\documentclass[a4paper,12pt]{article}

\usepackage[small,nohug,heads=vee]{diagrams}
\usepackage[latin1]{inputenc}
\usepackage[francais,english]{babel}

\usepackage[T1]{fontenc}
\usepackage[scaled]{helvet}

\usepackage{fullpage}
\usepackage{url}
\usepackage{a4,amsmath,amssymb,amsfonts,amsthm}
\usepackage{amsrefs}
\usepackage[retainorgcmds]{IEEEtrantools}

\newtheorem{theor}{Theorem}[section]
\newtheorem{prop}[theor]{Proposition}

\newtheorem{lemma}[theor]{Lemma}

\gdef\Td{\op{\rm Td}}
\gdef\ch{\op{\rm ch}}
\gdef\CB{{\mathcal B}}

\def\CF{{\mathcal F}}
\gdef\CE{{\mathcal E}}

\gdef\CL{{\mathcal L}}
\gdef\CO{{\mathcal O}}
\gdef\CP{{\mathcal P}}
\gdef\CS{{\mathcal S}}
\gdef\CM{{\mathcal M}}

\gdef\CA{{\mathcal A}}
\gdef\CH{{\mathcal H}}

\gdef\CT{{\mathcal T}}
\def\M#1{\mathbb#1}     
\def\mR{\M{R}}
\def\mZ{\M{Z}}

\def\mQ{\M{Q}}
\def\mC{\M{C}}

\gdef\beginProof{\par{\bf Proof. }}
\gdef\endProof{$\square$\par}
\gdef\ari#1{{\widehat{#1}}}
\gdef\mtr#1{\overline{#1}}
\gdef\c1{\op{{{\rm c}_1}}}
\gdef\ac1{\ari{\rm c}_1}
\gdef\aceq1{\ari{\rm c}_{{\rm eq},1}}

\gdef\ach{\ari{\rm ch}}
\gdef\aCH{\ari{\rm CH}}

\gdef\Spec{\op{\rm Spec}}

\gdef\refeq#1{(\ref{#1})}

\def\T{{\rm T}}

\def\R{{\rm R}}

\def\reg{{\rm reg}}

\def\m2{{\mu_2}}
\def\rk{{\rm rk}}

\def\op#1{\operatorname{#1}}

\def\be{\begin{equation}}
\def\ee{\end{equation}}

\def\wt#1{\widetilde{#1}}

\def\dR{{\rm dR}}

\def\an{{\rm an}}
\def\Id{{\rm Id}}
\def\log{{\rm log}}

\def\ddc{{\rm dd}^c}

\def\mfg{{\mathfrak g}}
\def\im{{\rm im}}

\def\cyc{{\rm cyc}}

\def\mod{{\rm mod}}
\def\BC{{\rm BC}}

\author{V. Maillot\footnote{Institut de Math\'ematiques de Jussieu,
Universit\'e Paris 7 Denis Diderot, C.N.R.S.,
Case Postale 7012,
2 place Jussieu,
F-75251 Paris Cedex 05, FRANCE,
E-mail : vmaillot@math.jussieu.fr}\ \ \&
D. R\"ossler\footnote{Institut de Mathématiques, \'Equipe Emile Picard,
Université Paul Sabatier,
118 Route de Narbonne,
31062 Toulouse cedex 9,
FRANCE, E-mail: rossler@math.univ-toulouse.fr}
}

\title{On a canonical class of Green currents for the unit sections of abelian schemes
\footnote{Mathematics Subject Classification (MSC2010): 14G40, 14K15, 11G16, 11G55, 14K25}}

\begin{document}

\maketitle

\begin{abstract}
We show that on any abelian scheme over a complex quasi-projective smooth variety, there
is a Green current for the zero-section, which is axiomatically determined up to $\partial$ and $\bar\partial$-exact differential
forms. On an elliptic curve, this current specialises to a Siegel function. We prove generalisations
of classical properties of Siegel functions, like distribution relations and reciprocity laws.
Furthermore, as an application of a refined version of the arithmetic Riemann-Roch theorem, we show that the above current, when restricted to a torsion section, is the realisation in analytic Deligne cohomology of an element of the (Quillen) $K_1$ group of the base, the corresponding denominator being given by the denominator of a Bernoulli number. This generalises the second Kronecker limit formula and the denominator
$12$ computed by Kubert, Lang and Robert in the case of Siegel units.  Finally, we prove an analog in Arakelov theory of a Chern class formula of Bloch and Beauville, where the canonical current plays a key role.
\end{abstract}

\bibliographystyle{plain}

\parindent=0pt

\begin{flushright}
\it For Christophe Soulé, on the occasion of his sixtieth birthday
\end{flushright}

\section{Introduction}

In this article, we show that on any abelian scheme over a complex quasi-projective smooth variety,
there is a Green current $\mfg$ for the the cycle given by the zero-section of
the abelian scheme, which is uniquely determined, up to $\partial$ and $\bar\partial$-exact
forms, by three axioms --- see Theorem \ref{mainth0}.

We proceed to show that this Green current is naturally compatible
with isogenies (ie it satisfies distribution relations), up to $\partial$ and $\bar\partial$-exact
forms (see Theorem \ref{mainth1}.1), and that it intervenes in an
Arakelov-theoretic generalization of a formula of Bloch and Beauville
(see \cite[p. 249]{Beauville-Sur-Chow} for the latter), which is proven here without resorting to the Fourier-Mukai transform. See Theorem
 \ref{mainth1}.2. Furthermore, we show that if the basis of the abelian scheme is a point (ie if the abelian scheme is an abelian variety), then
 $\mfg$ is a harmonic Green current (see Theorem \ref{mainth1}.3). The current
  $\mfg$ is also shown
 to be compatible with products (see Theorem \ref{mainth1}.5).

 Finally, we show that the restriction of $\mfg$ to the complement of the zero-section
 has a spectral interpretation. Up to a sign, it is given there by the degree $(g-1)$ part of the analytic
 torsion form of the Poincaré bundle of the abelian scheme. See Theorem \ref{mainth2}.1 for
 this. In point 2 of the same theorem, we show that the restriction of the higher analytic torsion form to
 torsion sections, which never meet the zero-section, lies in the rational image of the
 Beilinson regulator from $K_1$ to analytic Deligne cohomology and we give a multiplicative upper bound for the denominators involved.
 To prove Theorem \ref{mainth2}.1, we make heavy use of the arithmetic
 Riemann-Roch theorem in higher degrees proven in \cite{Gillet-Rossler-Soule-An-arithmetic} and to compute the denominators described
 in Theorem \ref{mainth2}.2, we apply the Adams-Riemann-Roch theorem
 in Arakelov geometry proven in \cite{Rossler-Adams}.

 If one specializes to elliptic schemes (i.e. abelian schemes of relative dimension $1$) the results proven in Theorems \ref{mainth0},
 \ref{mainth1} and \ref{mainth2} one recovers many results contained in the classical theory
 of elliptic units. In particular, on elliptic schemes the current $\mfg$
  is described by a Siegel function and the spectral interpretation
 of $\mfg$ specializes to the second Kronecker limit formula. The reciprocity law for elliptic units (ie the analytic description of the action of the Galois group on the elliptic units) is also easily obtained and (variants) of the results of Kubert-Lang and Robert on the fields of definition of elliptic units are recovered as a special case of the above denominator computations.   Details about elliptic schemes are given in section \ref{classun} where references to the classical literature are also given. The reader will notice that even in the case of elliptic schemes, our methods of proof are quite different from the classical ones.

The current $\mfg$ can also be used to describe the realisation in analytic Deligne cohomology of the degree $0$ part of the polylogarithm on abelian schemes introduced by J. Wildeshaus in \cite{Wildeshaus-Realizations} (see also \cite{Kings-K-theory}). The fact that this should be the case was a conjecture of G. Kings. His conjecture is proven in \cite{Kings-Rossler-Currents}.

 Here is a detailed description of the results.

Let $(R,\Sigma)$ be an arithmetic ring. By definition, this means
that $R$ is an excellent regular ring, which comes with a finite conjugation-invariant set $\Sigma$ of embeddings into $\mC$ (see \cite[3.1.2]{Gillet-Soule-Arithmetic}).
For example $R$ might be $\mZ$ with its unique embedding into $\mC$, or $\mC$ with
the identity and complex conjugation as embeddings.

Recall that an arithmetic variety over $R$ is a scheme, which is flat and of finite type over $R$.
In this text, all arithmetic varieties over $R$ will also be assumed to be regular, as well as
quasi-projective over $R$.
For any arithmetic variety $X$ over $R$, we write as usual
$$
X(\mC):=\coprod_{\sigma\in \Sigma}(X\times_{R,\sigma}\mC)(\mC)=\coprod_{\sigma\in \Sigma}X(\mC)_\sigma.
$$
Let
$D^{p,p}(X_\mR)$ (resp. $A^{p,p}(X_\mR)$) be the $\mR$-vector space of currents (resp. differential forms) $\gamma$ on $X(\mC)$ such that
\begin{itemize}
\item[$\bullet$] $\gamma$ is a real current (resp. differential form) of type $(p,p)$;
\item[$\bullet$] $F_\infty^*\gamma=(-1)^p\gamma$,
\end{itemize}
where $F_\infty:X(\mC)\to X(\mC)$ is the real analytic involution given by complex conjugation.
We then define
$$
\widetilde{D}^{p,p}(X_\mR):=D^{p,p}(X_\mR)/(\im\ \partial+\im\ \bar\partial)\
$$
(resp.
$$
\widetilde{A}^{p,p}(X_\mR):=A^{p,p}(X_\mR)/(\im\ \partial+\im\ \bar\partial)\ \qquad\quad).
$$

All these notations are standard in Arakelov geometry. See \cite{Soule-Lectures} for
a compendium. It is shown in \cite[Th. 1.2.2 (ii)]{Gillet-Soule-Arithmetic} that the natural map $\widetilde{A}^{p,p}(X_\mR)\to \widetilde{D}^{p,p}(X_\mR)$ is an injection.

If $Z$ a closed complex submanifold of $X(\mC)$, we shall write
more generally
$D^{p,p}_Z(X_\mR)$  for the $\mR$-vector space of currents  $\gamma$ on $X(\mC)$ such that
\begin{itemize}
\item[$\bullet$] $\gamma$ is a real current of type $(p,p)$;
\item[$\bullet$] $F_\infty^*\gamma=(-1)^p\gamma$;
\item[$\bullet$] the wave-front set of $\gamma$ is included in the real conormal bundle of $Z$ in $X(\mC)$.
\end{itemize}
Similarly, we then define the $\mR$-vector spaces
$$
\widetilde{D}^{p,p}_Z(X_\mR):=D^{p,p}_Z(X_\mR)/D^{p,p}_{Z,0}(X_\mR)
$$
where
$D^{p,p}_{Z,0}(X_\mR)$ is the set of currents $\gamma\in D^{p,p}_Z(X_\mR)$ such that: there exists
a complex current $\alpha$ of type $(p-1,p)$ and a complex current $\beta$ of type $(p,p-1)$  such that $\gamma:=\partial\alpha+\bar\partial\beta$ and such that the
wave-front sets of $\alpha$ and $\beta$ are included in the real conormal bundle of $Z$ in $X(\mC)$.

See \cite{Hoermander-The-analysis-I} for the definition (and theory) of the wave-front set.

It is a consequence of \cite[Cor. 4.7]{Burgos-Litcanu-Singular} that the natural morphism $\widetilde{D}^{p,p}_Z(X_\mR)\to \widetilde{D}^{p,p}(X_\mR)$ is an injection.\footnote{many thanks to J.-I. Burgos for bringing this to our attention} Thus the real vector space $\widetilde{D}^{p,p}_Z(X_\mR)$ can be
identified with a subspace of the real vector space $\widetilde{D}^{p,p}(X_\mR)$.

Furthermore, it is a consequence of \cite[Th. 4.3]{Burgos-Litcanu-Singular} that for any $R$-morphism $f:Y\to X$ of arithmetic varieties, there is a natural morphism
of $\mR$-vector spaces $$f^*:\widetilde{D}^{p,p}_Z(X_\mR)\to \widetilde{D}^{p,p}_{f(\mC)^*(Z)}(Y_\mR),$$ provided
$f(\mC)$ is transverse to $Z$. This morphism extends the morphism  $\widetilde{A}^{p,p}(X_\mR)\to \widetilde{A}^{p,p}(Y_\mR)$,
which is obtained by pulling back differential forms.

Fix now $S$ an arithmetic variety over
$R$.
Let $\pi:\CA\to S$ be an abelian scheme over $S$ of relative dimension $g$. We shall write
as usual $\CA^\vee\to S$ for the dual abelian scheme. Write $\epsilon$ (resp. $\epsilon^\vee$) for the zero-section of $\CA\to S$ (resp.
$\CA^\vee\to S$) and also $S_0$ (resp. $S_0^\vee$) for the image of $\epsilon$
(resp. $\epsilon^\vee$). We denote by the symbol $\CP$ the Poincaré bundle on $\CA\times_S\CA^\vee$.
We equip the Poincaré bundle $\CP$ with the unique
metric $h_\CP$ such that the canonical rigidification of $\CP$ along the
zero-section $\CA^\vee\to\CA\times_S\CA^\vee$ is an isometry and such that
the curvature form of $h_\CP$ is translation invariant along the fibres of the map
\mbox{$\CA(\mC)\times_{S(\mC)}\CA^\vee(\mC)\to\CA^\vee(\mC)$.} We write $\mtr{\CP}:=(\CP,h_\CP)$ for
the resulting hermitian line bundle. Write $\mtr{\CP}^0$ be the restriction of
$\CP$ to $\CA\times_S(\CA^\vee\backslash S_0^\vee)$.

The aim of this text is now to prove the following three theorems.

\begin{theor}
There is a unique class of currents $\mfg_\CA\in \widetilde{D}^{g-1,g-1}(\CA^\vee_\mR)$ with the
following three properties:
\begin{itemize}
\item[\rm (a)] Any element of $\mfg_\CA$ is a Green current for $S_0^\vee(\mC)$.
\item[\rm (b)] The identity
$
(S_0^\vee,\mfg_\CA)=(-1)^g p_{2,*}(\ari{\ch}(\mtr{\CP}))^{(g)}
$
 holds in $\ari{\rm CH}^g(\CA^\vee)_\mQ$.
\item[\rm (c)]  The identity $\mfg_\CA=[n]_*\mfg_\CA$ holds for all $n\geqslant 2$.
\end{itemize}\label{mainth0}
\end{theor}
Here the morphism $p_2$ is the second projection $\CA\times_S\CA^\vee\to\CA^\vee$ and
$[n]:\CA^\vee\to\CA^\vee$ is the multiplication-by-$n$ morphism.
The symbol $\ari{\ch}(\cdot)$ refers to the arithmetic Chern character and $\ari{\rm CH}^\bullet(\cdot)$ is the arithmetic
Chow group. See \cite[1.2]{Gillet-Soule-Arithmetic} for the notion of Green current.

{\bf Supplement.} The proof of Theorem \ref{mainth0} given below shows that if $S$ is assumed proper over $\Spec R$, then
the condition (b) can be replaced by the following weaker condition :

(b)' The identity of currents $\ddc\mfg+\delta_{S_0^\vee(\mC)}=(-1)^g p_{2,*}({\ch}(\mtr{\CP}))^{(g)}$
holds.

Here $\ddc:={i\over 2\pi}\partial\bar\partial$ and $\delta_{S_0^\vee(\mC)}$ is the Dirac current associated to
$S_0^\vee(\mC)$ in $\CA^\vee(\mC)$. Furthermore, ${\ch}(\mtr{\CP})$ is the Chern character form of
the hermitian bundle $\mtr{\CP}$. See \cite[Intro.]{Gillet-Soule-CharI} for this.

{\bf Remarks.} {\bf (1)} The condition (b) apparently makes the current $\mfg_\CA$ dependent on
the arithmetic structure of $\CA$. We shall show in \ref{mainth1}.4 below
that this is not the case. In particular, in defining $\mfg_\CA$, we could have assumed that
$R=\mC$. The settting of arithmetic varieties is used
in Theorem \ref{mainth0} because it is the natural one for property (b). Formula \refeq{dmfg0}
gives a purely analytic expression  for
$\mfg_\CA$. {\bf (2)} It is tempting to try to refine Theorem \ref{mainth0} by using in property (b) the arithmetic Chow groups defined by J.-I. Burgos (in \cite{Burgos-Arithmetic}) rather than the arithmetic Chow groups of Gillet-Soulé. One would then obtain a class of forms with certain logarithmic singularities, rather than a class of currents. Such a refinement does not seem to be easily attainable though, because of the lack of covariant functoriality of the spaces of forms mentioned in the last sentence.

The next theorem gives some properties of the class of currents $\mfg_\CA$.

Let ${\CL}$ be a rigidified line bundle on $\CA$. Endow
$\CL$ with the unique hermitian metric $h_\CL$, which is compatible with the rigidification and whose curvature form is translation-invariant on the fibres
of $\CA(\mC)\to S(\mC)$. Let $\mtr{\CL}:=(\CL,h_\CL)$ be the resulting hermitian line bundle. Let
$\phi_\CL:\CA\to\CA^\vee$ be the polarisation morphism induced by $\CL$.
\begin{theor}
\begin{enumerate}
\item Let $\iota:\CA\to\CB$ be an isogeny of abelian schemes over $S$. Then  the identity $\iota^\vee_*(\mfg_\CB)=\mfg_\CA$ holds.
\item  Suppose that $\CL$ is ample relatively to $S$ and symmetric.  Then the equalities
$$
(S_0^\vee,\mfg_\CA)=(-1)^g p_{2,*}(\ari{\ch}(\mtr{\CP}))={1\over g!\sqrt{\deg(\phi_\CL)}}\phi_{\CL,*}(\ac1(\mtr{\CL})^g)
$$
are verified in $\ari{\rm CH}^\ast(\CA^\vee)_{\mQ}$.
\item If $S\to\Spec R$ is the identity on $\Spec R$  then any element of $\mfg_\CA$ is a harmonic Green current for $S_0^\vee(\mC)$,  where $\CA^\vee(\mC)$ is endowed with a conjugation invariant Kähler metric, whose Kähler form is translation invariant.
\item The class $\mfg_\CA$ is invariant under any change of arithmetic rings $(R,\Sigma)\to(R',\Sigma')$.
\item Let $\CB\to S$ be another abelian scheme and let $\pi_{\CA^\vee}:\CA^\vee\times_S\CB^\vee\to \CA^\vee$
(resp.  $\pi_{\CB^\vee}:\CA^\vee\times_S\CB^\vee\to \CB^\vee$) be the natural projections. Then
\begin{equation}
\mfg_{\CA\times_S\CB}=\pi_{\CA^\vee}^*(\mfg_{\CA})\ast \pi_{\CB^\vee}^*(\mfg_{\CB})
\label{eqast}
\end{equation}
\item The class of currents $\mfg_\CA$ lies in $\widetilde{D}^{g-1,g-1}_{S_0^\vee(\mC)}(\CA^\vee_\mR)$.
\item Let $T$ be a an arithmetic variety over $R$ and let $T\to S$ be a morphism of schemes over $R$. Let
$\CA_T$ be the abelian scheme obtained by base-change and let $\BC:\CA_T\to \CA$ be the corresponding
morphism. Then $\BC(\mC)$ is tranverse to $S_0^\vee(\mC)$ and $\BC^*\mfg_\CA=\mfg_{\CA_T}$.
\end{enumerate}
\label{mainth1}
\end{theor}
Here $\iota^\vee:\CB^\vee\to\CA^\vee$ is the isogeny, which is dual to $\iota$.
For the notion of harmonic Green current, see \cite{Bost-Green} and \cite{Kunnemann-Arakelov}. The pairing $\ast$ appearing in the equation \refeq{eqast}  is the $\ast$-product of Green currents. See \cite[par. 2.2.11, p. 122]{Gillet-Soule-Arithmetic} for the definition.

Recall that an $S$-isogeny between the abelian schemes $\CA$ and $\CB$ is a flat and finite $S$-morphism
$\CA\to \CB$, which is compatible with the group-scheme structures. The symbol
$\ac1(\cdot)$ refers the first arithmetic Chern class; see \cite[Intro.]{Gillet-Soule-CharI} for this notion.

Theorem \ref{mainth1}.1 generalizes to higher degrees the distribution
relations of Siegel units. See section \ref{classun}  below for details. If the morphism
$S\to\Spec R$ is the identity on $\Spec R$ and $R$ is the ring of integers of a number field, then
it is shown in \cite[Prop. 11.1 (ii)]{Kunnemann-Arakelov} that Theorem \ref{mainth1}.3
implies Theorem \ref{mainth1}.2. Still in the situation where $S\to\Spec R$ is the identity on $\Spec R$, another
construction of a Green current for $S_0^\vee(\mC)$ is described in A. Berthomieu's thesis \cite{Berthomieu-These}. The current constructed by Berthomieu is likely to be harmonic
(it is not proven in \cite{Berthomieu-These}, but according to the author [private communication] it can easily be shown).  The current constructed in \cite{Berthomieu-These} satisfies the identity in Theorem \ref{mainth2}.1 by construction.

The last theorem relates the current $\mfg_\CA$ to the Bismut-Köhler analytic torsion form of
the Poincaré bundle (see \cite[Def. 3.8, p. 668]{Bismut-Koehler-Higher} for the definition).

Let $\lambda$ be a $(1,1)$-form on $\CA(\mC)$ defining a Kähler fibration structure on
the fibration $\CA(\mC)\to S(\mC)$ (see \cite[par. 1]{Bismut-Koehler-Higher} for this notion). With the form
$\lambda$, one can canonically associate a hermitian metric on the relative cotangent bundle $\Omega_{\CA/S}$ and we shall write
$\mtr{\Omega}_{\CA/S}$ for the resulting hermitian vector bundle.
We suppose that $\lambda$ is translation invariant on the fibres of the map $\CA(\mC)\to S(\mC)$ as well as conjugation invariant.
 We shall write
 $$T(\lambda,\mtr{\CP}^0)\in\widetilde{A}((\CA^\vee\backslash S_0^\vee)_\mR):=
 \bigoplus_{p\geqslant 0}\widetilde{A}^{p,p}((\CA^\vee\backslash S_0^\vee)_\mR)$$ for
the Bismut-Köhler higher analytic torsion form of $\mtr{\CP}^0$ along the fibration
$$\CA(\mC)\times_{S(\mC)}(\CA^\vee(\mC)\backslash S_0^\vee(\mC))\longrightarrow\CA^\vee(\mC)\backslash S_0^\vee(\mC).$$
For any regular arithmetic variety $X$ over $R$, the (Beilinson) regulator map
gives rise to a morphism of groups
$$
\reg_\an:K_1(X)\longrightarrow\bigoplus_{p\geqslant 0}H^{2p-1}_{D,\an}(X_\mR,\mR(p)).
$$
To define the space $H^{2p-1}_{D,\an}(X_\mR,\mR(p))$ and the map $\reg_\an$, let us first
write $H^{*}_{D,\an}(X,\mR(\cdot))$ for the analytic real Deligne cohomology of $X(\mC)$.
By definition,
$$
H^{q}_{D,\an}(X,\mR(p)) := \M{H}^q(X(\mC),\mR(p)_{D,\an})
$$
where $\mR(p)_{D,\an}$ is the complex of sheaves of $\mR$-vector spaces
$$
0\to\mR(p)\to{\cal O}_{X(\mC)}\stackrel{d\;}{\to}\Omega^1_{X(\mC)}\to\dots\to\Omega^{p-1}_{X(\mC)}\to 0
$$
on $X(\mC)$ (for the ordinary topology). Here $\mR(p):=(2i\pi)^p\,\mR\subseteq\mC$. We now define
$$
H^{2p-1}_{D,\an}(X_\mR,\mR(p)):=\{\gamma\in H^{2p-1}_{D,\an}(X,\mR(p))\;|\;F_\infty^*\gamma=(-1)^p\gamma\}.
$$
By construction, the regulator map $K_1(X)\to\oplus_{p\geqslant 0}H^{2p-1}_{D,\an}(X,\mR(p))$
(see \cite{Burgos-Wang-Higher} for a direct construction of the regulator and further
references) factors through $\oplus_{p\geqslant 0}H^{2p-1}_{D,\an}(X_\mR,\mR(p))$
and thus gives rise to a map $K_1(X)\to\oplus_{p\geqslant 0}H^{2p-1}_{D,\an}(X_\mR,\mR(p))$.
This is the definition of the map $\reg_\an$.

It is shown in \cite[par 6.1]{Burgos-Wang-Higher} that there is a natural inclusion
\mbox{$H^{2p-1}_{D,\an}(X_\mR,\mR(p))\hookrightarrow \widetilde{A}^{p,p}(X_\mR)$.}

For the next theorem, define
$$N_{2g}:=2\cdot{\rm denominator}\,[(-1)^{g+1}B_{2g}/(2g)],$$ where $B_{2g}$ is
the $2g$-th Bernoulli number. Recall that the Bernoulli numbers are defined by the identity of power series:
$$
\sum_{t\geqslant 1}B_t{u^t\over t!}={u\over \exp(u)-1}.
$$
\begin{theor}
\begin{enumerate}
\item The class of differential forms $\Td(\epsilon^*\mtr{\Omega}_{\CA/S})\cdot T(\lambda,\mtr{\CP}^0)$ lies
in $\widetilde{A}^{g-1,g-1}((\CA^\vee\backslash S_0^\vee)_\mR)$ and the equality $$\mfg_\CA|_{\CA^\vee(\mC)\backslash S_0^\vee(\mC)}=
(-1)^{g+1}\Td(\epsilon^*\mtr{\Omega}_{\CA/S})\cdot T(\lambda,\mtr{\CP}^0)$$ holds. In particular
$T(\lambda,\mtr{\CP}^0)^{(g-1)}$ does not depend on $\lambda$.
\item Suppose that $\lambda$ is the first Chern form of a relatively ample rigidified line
bundle, endowed with its canonical metric.
Let \mbox{$\sigma\in\CA^\vee(S)$} be an element of finite order $n$, such that $\sigma^*S_0^\vee=\emptyset$.
Then
\begin{equation*}
g\cdot n\cdot N_{2g}\cdot \sigma^*T(\lambda,\mtr{\CP}^0)\in\op{image}(\reg_\an(K_1(S))).
\end{equation*}
\end{enumerate}
\label{mainth2}
\end{theor}
A the end of section \ref{proof132} (see the end of the proof of Lemma \ref{comlem}), we give a statement, which is slightly stronger than
Theorem \ref{mainth2}.2 (but more difficult to formulate).

Theorem \ref{mainth2}.1 can be viewed as a  generalization to higher degrees of
the second Kronecker limit formula (see \cite[chap. 20, par. 5, p. 276]{Lang-Elliptic} for the latter). Theorem \ref{mainth2}.2 generalizes to higher degrees part of
a classical statement on elliptic units and their fields of definition. See section \ref{classun} below.

{\bf Remark.} It would be interesting to have an analogue of Theorem \ref{mainth2}.2,
where $\reg_\an$ is replaced by the analytic cycle class $\cyc_\an$ (see \refeq{arfundseq} below). If
$S\simeq\Spec R$ and $R$ is  the ring of integers in a number field, then $\reg_\an$ and
$\cyc_\an$ can be identified but this is not true in general. In particular, this suggests that
the Bernoulli number  $12=N_2/2$, which appears in the denominators of elliptic units (see the last section) should be understood as coming from the natural integral
structure of the group $K_1(\cdot)_\mQ$ and not from the natural integral structure of the corresponding motivic cohomology
group $\bigoplus_{p}H^{2p-1}_{\cal M}(\cdot,\mZ(p))_\mQ$.


Some of the results of this article were announced in \cite{Maillot-Rossler-Elements}.

We shall provide many bibliographical references to ease the reading but the reader of this article is nevertheless assumed to have some familiarity with the language of
Arakelov theory, as expounded for instance in \cite{Soule-Lectures}.

{\bf Acknowledgments.} We thank J.-I. Burgos for patiently listening to our explanations
on the contents of this article and for his (very) useful comments over a number of years. We also thank C. Soulé, G. Kings, as well as S. Bloch and A. Beilinson for their interest. J. Kramer made several interesting remarks on the contents of this article and his input was very useful. Many thanks also to K. Köhler for his explanations on the
higher analytic torsion forms of abelian schemes. Finally, we are grateful to J. Wildeshaus for his feedback and for answering many questions on abelian polylogarithms.

{\bf Notations.} Here are the main notational conventions. Some of them
have already been introduced above. Recall that we wrote $\pi:\CA\to S$ for the structure morphism of
the abelian scheme $\CA$ over $S$. We also write $\pi^\vee:\CA^\vee\to S$ for the structure morphism of
the abelian scheme $\CA^\vee$ over $S$. Write $\mu=\mu_\CA:\CA\times_S\CA\to \CA$ for the addition morphism
and $p_1:\CA\times\CA^\vee\to\CA$,  $p_2:\CA\times\CA^\vee\to\CA^\vee$ for the obvious projections. We shall also also write
${\underbar p}_1,{\underbar p}_2:\CA\times\CA\to\CA$ and ${\underbar p}^\vee_1,{\underbar p}^\vee_2:\CA^\vee\times\CA^\vee\to\CA^\vee$ for more obvious projections. Recall that we wrote $\epsilon$ (resp. $\epsilon^\vee$) for the zero-section of $\CA\to S$ (resp.
$\CA^\vee\to S$) and also that we wrote $S_0$ (resp. $S_0^\vee$) for the image of $\epsilon$
(resp. $\epsilon^\vee$). Write $\omega_\CA:=\det(\Omega_{\CA/S})$ for the determinant
of the sheaf of differentials of $\CA$ over $S$. We let  $\ddc:={i\over 2\pi}\partial\bar\partial$.

\section{Proof of Theorem \ref{mainth0}}
\label{pmth0}

If $M$ is a smooth complex quasi-projective variety, we shall write $H^{*}_{D}(M,\mR(\cdot))$ for the
Deligne-Beilinson cohomology of $M$. We recall its definition. Let $\bar M$ be a smooth complex
projective variety, which contains $M$ as an open subscheme. We call
$\bar M$ a compactfication of $M$. Suppose furthermore
that $\bar M\backslash M$ is the underlying set of a reduced divisor with normal crossings $D$.
From now on, we view $M$ and $\bar M$ as complex analytic spaces and we work in
the category of complex analytic spaces.
Let $j:M\hookrightarrow \bar M$ be the given open embedding. There is a natural subcomplex $\Omega^\bullet_{\bar M}(\log\ D)$ of $j_*\Omega^\bullet_M$, called the {\it complex of holomorphic differential forms on $M$ with
logarithmic singularities along $D$.} The objects of $\Omega^\bullet_{\bar M}(\log\ D)$ are
locally free sheaves. We redirect the reader to  \cite[chap. 10]{Burgos-The-regulators} for the definition and further bibliographical references. Write
$F^p\Omega^\bullet_{\bar M}(\log\ D)$ for the subcomplex
$$
\Omega^p_{\bar M}(\log\ D)\longrightarrow \Omega^{p+1}_{\bar M}(\log\ D)\longrightarrow\cdots
$$
of $\Omega^\bullet_{\bar M}(\log\ D)$. Write $f^0_p:F^p\Omega^\bullet_{\bar M}(\log\ D)\to
j_*\Omega^\bullet_M$ for the inclusion morphism. Abusing notation, we shall identify
$\R j_*\mR(p)$ with the complex, which is the image by $j_*$ of the canonical flasque resolution of $\mR(p)$.
Similarly, we shall write $\R j_*\Omega^\bullet_M$ for the simple complex
associated to the image by $j_*$ of the canonical flasque resolution of $\Omega^\bullet_M$ (
the latter being  a double complex).
Write $f_p:F^p\Omega^\bullet_{\bar M}(\log\ D)\to
\R j_*\Omega^\bullet_M$ for the morphism obtained by composing $f^0_p$ with
the canonical morphism $j_*\Omega^\bullet_M\to \R j_*\Omega^\bullet_M$.
There is a natural morphism of complexes $\mR(p)\to\Omega^\bullet_M$
(where $\mR(p)$ is viewed as a complex with one object sitting in degree $0$) and by the functoriality
of the flasque resolution, we obtain a morphism $r_p:\R j_*\mR(p)\to\R j_*\Omega^\bullet_M$.
We now define the complex
$$
\mR(p)_D:={\rm simple}\big(\R j_*\mR(p)\oplus F^p\Omega^\bullet_{\bar M}(\log\ D)\stackrel{u_p}{\longrightarrow}
\R j_*\Omega_M^\bullet\big)
$$
where $u_p:=r_p-f_p$. By definition, Deligne-Beilinson cohomology is now defined by the
formula
$$
H^q_D(M,\mR(p)):=\M{H}^q(M,\mR(p)_D).
$$
Notice that by construction, $H^q_D(M,\mR(p))=H^q_{D,\an}(M,\mR(p))$ if
$M$ is compact (so that $D$ is empty). More generally, there is a natural "forgetful" morphism
of $\mR$-vector spaces $H^q_D(M,\mR(p))\to H^q_{D,\an}(M,\mR(p))$ (what is forgotten is the
logarithmic structure) ; see \cite[before Prop. 1.3]{Burgos-Arithmetic} for this.
It can be proven that Deligne-Beilinson cohomology does not depend on the choice of
the compactification $\bar M$. By its very definition, we have a canonical long exact sequence of $\mR$-vector spaces
\begin{equation}
\cdots\to H^{q-1}(M,\mC)\to H^q_D(M,\mR(p))\to H^q(M,\mR(p))\oplus F^p H^q(M,\mC)\to\cdots
\label{delong}
\end{equation}
where $F^p H^q(M,\mC)$ is the $p$-th term of the Hodge filtration of the mixed Hodge structure on
$H^q(M,\mC)$.
Furthermore the $\mR$-vector space $H^q_D(M,\mR(p))$
has a natural structure of contravariant functor from the category of
smooth quasi-projective varieties over $\mC$ to the category of $\mR$-vector spaces (see
\cite[Prop. 1.3]{Burgos-Arithmetic} for this).  If we equip the singular cohomology spaces $H^q(M,\mC)$ and
$H^q_D(M,\mR(p))$ with their natural contravariant structure, then the sequence \refeq{delong} becomes
an exact sequence of functors.

If $X$ is an arithmetic variety, then we define
$$
H^q_D(X_\mR,\mR(p)):=\{\gamma\in H^q_D(X(\mC),\mR(p))\;|\;F^{*}_{\infty}\gamma=(-1)^p\gamma\}.
$$

Before beginning with the proof of Theorem \ref{mainth0}, we shall prove the following key lemma.

\begin{lemma}
Let $n\geqslant 2$.
The eigenvalues of the $\mR$-endomorphism $[n]^*$ of the Deligne-Beilinson cohomology $\mR$-vector space $H^{2p-1}_{D}(\CA^\vee(\mC),\mR(p))$
lie in the set $\{1,n,n^2,\dots,n^{2p-1}\}$.
\label{keylem}
\end{lemma}
\beginProof
The existence of the exact sequence of functors \refeq{delong} shows that we have the following
exact sequence of $\mR$-vector spaces
$$
H^{2p-2}(\CA^\vee(\mC),\mC)\to H^{2p-1}_D(\CA^\vee(\mC),\mR(p))\to H^{2p-1}(\CA^\vee(\mC),\mR(p))\oplus
F^{p}H^{2p-1}(\CA^\vee(\mC),\mC)
$$
and that the differentials in this sequence are compatible with the natural contravariant action
of $[n]$.
Hence it is sufficient to prove the conclusion of the lemma for the $\mR$-vector spaces  $H^{2p-2}(\CA^\vee(\mC),\mC)$, $H^{2p-1}(\CA^\vee(\mC),\mR(p))$ and $H^{2p-1}(\CA^\vee(\mC),\mC)$. These spaces
are more easily tractable and can be approximated by the Leray spectral sequence
$$
E_2^{rs}=H^{r}(S(\mC),\R^s\pi^\vee(\mC)_*(K))\Rightarrow H^{r+s}(\CA^\vee(\mC),K)
$$
where $K$ is $\mR$ or $\mC$. Now notice that since $[n]$ is an $S$-morphism, this
spectral sequence carries a natural contravariant action of $[n]$, which is compatible with
the aforementionned action of $[n]$ on its abutment.
Consider the index $2p-1$.
We know that $[n]^*$ acts on $\R^s\pi^\vee(\mC)_*(K)$ by multiplication
by $n^s$. This may be deduced from known results on abelian varieties using the proper base-change theorem. Hence $[n]^*$ acts on $H^{r}(S(\mC),\R^s\pi^\vee(\mC)_*(K))$ by multiplication
by $n^s$ as well. Now the existence of the spectral sequence
shows that $H^{2p-1}(\CA^\vee(\mC),K)$ has
a natural filtration, which consists of subquotients of the spaces $H^{r}(S(\mC),\R^s\pi^\vee(\mC)_*(K))$, where
$r+s=2p-1$. Since $s\leqslant 2p-1$, this proves the assertion for the index $2p-1$.
The index $2p-2$ can be treated in an analogous fashion.
\endProof

{\bf Proof of uniqueness.}
Let $\mfg_\CA$ and $\mfg_\CA^0$ be elements of $\widetilde{D}^{g-1,g-1}(\CA^\vee_\mR)$ satisfying (a), (b) and (c).
Let $\kappa_\CA:=\mfg_\CA^0-\mfg_\CA\in\widetilde{D}^{g-1,g-1}(\CA^\vee)$ be the error term.

Recall the fundamental exact sequence
\begin{equation}
{\rm CH}^{g,g-1}(\CA^\vee)
\xrightarrow{\cyc_\an}\widetilde{A}^{g-1,g-1}(\CA^\vee_\mR)\stackrel{a\;}{\to}\ari{\rm CH}^g(\CA^\vee)\to
{\rm CH}^g(\CA^\vee)\to 0
\label{arfundseq}
\end{equation}
(see \cite[th. 3.3.5]{Gillet-Soule-Arithmetic} for this). Here ${\rm CH}^{g,g-1}(\cdot)$ is Gillet-Soul\'e's version of
one of Bloch's higher Chow groups. The group $\ari{\rm CH}^g(\CA^\vee)$ is the $g$-th arithmetic Chow group
and ${\rm CH}^g(\CA^\vee)$ is the $g$-th ordinary Chow group. By construction, there are maps
$$
{\rm CH}^{g,g-1}(\CA^\vee)\xrightarrow{\cyc}H^{2g-1}_D(\CA^\vee_\mR,\mR(g))\xrightarrow{\rm forgetful~} H^{2g-1}_{D,\an}(\CA^\vee_\mR,\mR(g)){\to}\widetilde{A}^{g-1,g-1}(\CA^\vee_\mR)
$$
whose composition is $\cyc_\an$. Here $\cyc$ is the cycle class map into Deligne-Beilinson cohomology;
the second map from the left is the forgetful map and the third one is the natural inclusion mentioned
before Theorem \ref{mainth2}.

Now let $n\geqslant 2$. Let
\begin{equation}
V:={\rm image}\big(H^{2g-1}_D(\CA^\vee_\mR,\mR(g))\xrightarrow{\rm forgetful~} H^{2g-1}_{D,\an}(\CA^\vee_\mR,\mR(g))\big)
\label{defV}
\end{equation}
be the image of the forgetful map from $H^{2g-1}_D(\CA^\vee_\mR,\mR(g))$ to $H^{2g-1}_{D,\an}(\CA^\vee_\mR,\mR(g))$. By (b), we know that $\kappa_\CA\in V$.  Furthermore
$V$ is invariant under $[n]^*$. In fact, by Lemma \ref{keylem}, $[n]^*$ restrict to an injective
morphism $V\to V$, which is thus an isomorphism, since $V$ is finite dimensional.
Now the projection formula shows that the equation $[n]_*[n]^*=n^{2g}$ is valid in $H^{2g-1}_{D,\an}(\CA^\vee_\mR,\mR(g))$ and we conclude that that $V$ is also invariant under $[n]_*$. The same equation $[n]_*[n]^*=n^{2g}$ now shows that
the eigenvalues of $[n]_*$ on $V$ lie in the set
$\{n^{2g},n^{2g-1},\dots,n\}$. In particular, $[n]_*$ has no non-vanishing fixed point in $V$. Since
$[n]_*\kappa_\CA=\kappa_\CA$ by (c), this
proves that $\kappa_\CA=0$.

{\bf Proof of existence.} As very often, the proof of existence is inspired by the proof of uniqueness.
Let $\mfg\in \widetilde{D}^{g-1,g-1}(\CA^\vee_\mR)$ be a class of Green currents for $S_0^\vee$ satisfying (b). To see that there is such a $\mfg$,
pick any Green current $\mfg'$ for $S_0^\vee(\mC)$, such that $F_\infty^*\mfg'=(-1)^{g-1}\mfg'$.
This exists by
\cite[th.1.3.5]{Gillet-Soule-Arithmetic}. Now a basic property of the Fourier-Mukai transformation for abelian
schemes (see \cite[Lemme 1.2.5]{Laumon-Transformation}) implies that $(-1)^g p_{2,*}({\ch}({\CP}))^{(g)}=S_0^\vee$ in ${\rm CH}^g(\CA^\vee)_\mQ$.
Hence, looking at the sequence \refeq{arfundseq}, we see that
there exists $\alpha\in\widetilde{A}^{g-1,g-1}(\CA^\vee_\mR)$ such that $(a\otimes\mQ)(\alpha)=(S_0^\vee,\mfg')-
(-1)^g p_{2,*}(\ari{\ch}(\mtr{\CP}))^{(g)}$. If we define $\mfg:=\mfg'-\alpha$, we obtain the required class of Green currents.
Now fix $n\geqslant 2$ and let $c:=\mfg-[n]_*\mfg$. We shall prove below (see \refeq{impeq1})
that
$$
[n]_*p_{2*}(\ari{\ch}(\mtr{\CP}))^{(g)}=p_{2*}(\ari{\ch}(\mtr{\CP}))^{(g)}.
$$
This implies that $c$ lies in the space $V$ defined in \refeq{defV} above.
Now recall that we proved that $[n]_*$ sends $V$ on $V$ and that $1$ is not a root of the characteristic
polynomial of $[n]_*$ as an endomorphism of $V$. Hence the linear equation in $x$
$$
x-[n]_*x=c
$$
has a unique solution in $V$. Call this solution $c_0$. By construction the current
$\mfg_0:=\mfg+c_0$ satisfies the equation $\mfg_0-[n]_*\mfg_0=0$. Now let $m\geqslant 2$ be
another natural number. We have seen above that $\mfg_0-[m]_*\mfg_0$ lies in $V$. On the other
hand
$$
[n]_*(\mfg_0-[m]_*\mfg_0)=[n]_*\mfg_0-[m]_*[n]_*\mfg_0=\mfg_0-[m]_*\mfg_0
$$
hence $\mfg_0-[m]_*\mfg_0$ is a fixed point of $[n]_*$ in $V$. This implies that $\mfg_0-[m]_*\mfg_0=0$.
This proves that $\mfg_0$ satisfies (a), (b) and (c).

\section{Proof of Theorem \ref{mainth1}}

\subsection{Proof of \ref{mainth1}.1}

By the definition of the dual isogeny, there is a diagram
\begin{diagram}
\CA\times_S\CB^\vee &\rTo^{\;\iota\times\Id\;\;} & \CB\times_S \CB^\vee & \rTo & \CB^\vee\\
\dTo^{\Id\times\iota^\vee}                          &        &                              &         & \dTo^{\iota^\vee}\\
\CA\times_S\CA^\vee &        &\rTo                       &         & \CA^\vee
\end{diagram}
such that
$$
(\Id\times\iota^\vee)^*\CP_\CA\simeq (\iota\times\Id)^*\CP_\CB
$$
and such that the outer square is cartesian. Here $\CP_\CA:=\CP$ and $\CP_\CB$ is the Poincaré bundle
of $\CB$ over $S$.

First notice that the arithmetic Riemann-Roch theorem \cite{Gillet-Rossler-Soule-An-arithmetic}
implies that
$$
\ach((\iota\times\Id)_*(\mtr{\CO}_\CA))=\deg\ \iota
$$
in $\ari{\rm CH}^\bullet(\CB\times_S\CB^\vee)_\mQ$.

Now we compute
\begin{eqnarray}
\iota^{\vee,*} p_{2*}^\CA(\ari{\ch}(\mtr{\CP}_\CA))=
p_{2*}^\CB(\ari{\ch}(\mtr{\CP}_\CB)\ari{\ch}((\iota\times\Id)_*(\mtr{\CO}_\CA)))=
(\deg\ \iota)\cdot p_{2*}^\CB(\ari{\ch}(\mtr{\CP}_\CB)).
\label{ioteq}
\end{eqnarray}
Here we used the projection formula for arithmetic Chow theory (see
\cite[Th. 4.4.3, 7.]{Gillet-Soule-Arithmetic}) and the
fact that the push-forward map in arithmetic Chow theory commutes with base-change.
We may now compute
\begin{eqnarray*}
\iota^\vee_*\iota^{\vee,*} p_{2*}^\CA(\ari{\ch}(\mtr{\CP}_\CA))=
(\deg\ \iota)\cdot p_{2*}^\CA(\ari{\ch}(\mtr{\CP}_\CA))=
(\deg\ \iota)\cdot \iota^\vee_*p_{2*}^\CB(\ari{\ch}(\mtr{\CP}_\CB))).
\end{eqnarray*}
In other words, we have
\begin{equation}
p_{2*}^\CA(\ari{\ch}(\mtr{\CP}_\CA))=
\iota^\vee_*p_{2*}^\CB({\ari{\ch}}(\mtr{\CP}_\CB))).
\label{impeq1}
\end{equation}
Furthermore, since $\iota^\vee$ restricts to an isomorphism between the zero-sections,
the class of currents $\iota_*^\vee(\mfg_\CB)$ consists of Green currents
for $S_{0,\CA}^\vee$. All this shows that $\mfg_\CA-\iota_*^\vee(\mfg_\CB)$ lies inside
the space $V$ defined in \ref{defV}. Recall that $V$ is the image of the forgetful map
$H^{2g-1}_D(\CA^\vee_\mR,\mR(g))\xrightarrow{\rm forgetful~} H^{2g-1}_{D,\an}(\CA^\vee_\mR,\mR(g))$.
To conclude, notice that for any $n\geqslant 2$, we have
$$
[n]_*(\mfg_\CA-\iota_*^\vee(\mfg_\CB))=\mfg_\CA-\iota_*^\vee(\mfg_\CB)
$$
since $[n]$ commutes with $\iota^\vee$. It was shown just before the proof of existence in the proof of Theorem
\ref{mainth0} that $[n]_*$ leaves $V$ invariant and has no non-vanishing fixed points in $V$.
Thus $\mfg_\CA-\iota_*^\vee(\mfg_\CB)=0$.

\subsection{Proof of \ref{mainth1}.2}

We shall prove the equivalent identities
\begin{equation}
{(-1)^g\over g!\sqrt{\deg(\phi_\CL)}}\phi_{\CL,*}(\ac1(\mtr{\CL})^g)=p_{2,*}(\ari{\ch}(\mtr{\CP}))^{(g)}
\label{poincample}
\end{equation}
and
\begin{equation}
p_{2,*}(\ari{\ch}(\mtr{\CP}))^{(k)}=0
\label{poincample1}
\end{equation}
if $k\not= g$.

For the equality \refeq{poincample1}, notice that in view of \refeq{ioteq} and the fact that
$[n]^\vee=[n]$, we have
\begin{equation}
[n]^*(p_{2,*}(\ari{\ch}(\mtr{\CP})))=n^{2g}\cdot p_{2,*}(\ari{\ch}(\mtr{\CP}))
\label{eqho1}
\end{equation}
for any $n\geqslant 2$. On the other hand, since $(\Id\times [n])^*\bar\CP=\bar \CP^{\otimes n}$, we have
also
\begin{equation}
[n]^*(p_{2,*}(\ari{\ch}(\mtr{\CP})))=\sum_{k\geqslant 0}n^{k+g}\cdot p_{2,*}(\ari{\ch}(\mtr{\CP}))^{(k)}
\label{eqho2}
\end{equation}
and comparing equations \refeq{eqho1} and \refeq{eqho2} as polynomials in $n$ proves
equation \refeq{poincample1}.

We now proceed to the proof of equation \refeq{poincample}.
Notice that the line bundle $\mu^*{\CL}\otimes \underbar p_1^*{\CL}^\vee\otimes \underbar p_2^*{\CL}^\vee$ on $\CA\times_S\CA$ carries a natural rigidification on the
zero section
$\CA\xrightarrow{(\Id,\epsilon)}\CA\times_S\CA$ and that the same line bundle is algebraically equivalent to $0$ on each
geometric fibre of the morphism $\underbar p_2:\CA\times_S\CA\to\CA$. Hence there is a unique morphism
$\phi_\CL:\CA\to\CA^\vee$, the polarisation morphism induced by $\CL$, such that there is an isomorphism
of rigidified line bundles
\begin{equation}
(\Id\times\phi_\CL)^*{\CP}\simeq\mu^*{\CL}\otimes \underbar p_1^*{\CL}^\vee\otimes \underbar p_2^*{\CL}^\vee .
\label{eqpol}
\end{equation}
Furthermore, if we endow the line bundles on both sides of \refeq{eqpol} with their natural metrics, this
isomorphism becomes an isometry, because both line bundles carry metrics that are compatible with the rigidification and the curvature forms of both sides are translation invariant (in fact $0$) on the fibres
of the map $\underbar p_2(\mC)$.

We shall now give a more concrete expression for $\underbar p_{2*}(\ari{\ch}(\mu^*\mtr{\CL})\ari{\ch}(\underbar p_1^{*}\mtr{\CL}^\vee)\ari{\ch}(\underbar p_2^*\mtr{\CL}^\vee))$. We first
make the calculation
\begin{IEEEeqnarray*}{rCl}
&&\underbar p_{2*}(\ari{\ch}(\mu^*\mtr{\CL})\ari{\ch}(\underbar p_1^{*}\mtr{\CL}^\vee)\ari{\ch}(\underbar p_1^*\mtr{\CL}))=
\underbar p_{2*}(\ari{\ch}(\mu^*\mtr{\CL}))=
\underbar p_{2*}(\alpha^*\underbar p_1^{*}\ari{\ch}(\mtr{\CL}))=
\underbar p_{2*}(\underbar p_1^{*}\ari{\ch}(\mtr{\CL}))
\end{IEEEeqnarray*}
where $\alpha:\CA\times_S\CA\to\CA\times_S\CA$ is the $\underbar p_2$-automorphism
$\alpha:=(\mu,\underbar p_2)$. Now notice that for any $n\geqslant 2$,
\begin{eqnarray*}
&&\underbar p_{2*}(([n]\times\Id)_*([n]\times\Id)^*\underbar p_1^{*}\ari{\ch}(\mtr{\CL}))=
n^{2g}\underbar p_{2*}(\underbar p_1^{*}\ari{\ch}(\mtr{\CL}))\\&=&
\underbar p_{2*}(([n]\times\Id)^*\underbar \underbar \underbar p_1^{*}\ari{\ch}(\mtr{\CL}))=
\underbar p_{2*}(\sum_{l\geqslant 1}n^{2l}(\underbar p_1\ari{\ch}(\mtr{\CL}))^{(l)}).
\end{eqnarray*}
Here we used the isometric isomorphism $[n]^*\mtr{\CL}\simeq\mtr{\CL}^{\otimes n^2}$ (recall
that $\CL$ is symmetric). We deduce that
$$
\underbar p_{2*}(\underbar p_1^{*}\ari{\ch}(\mtr{\CL}))=\underbar p_{2*}(\underbar p_1^{*}\ari{\ch}(\mtr{\CL})^{(g)})
=\sqrt{\deg(\phi_{\CL})}
$$
Thus, using the projection formula, we see that
$$
\underbar p_{2*}(\ari{\ch}(\underbar p_1^{*}\mtr{\CL})\ari{\ch}(\mu^*\mtr{\CL})\ari{\ch}(\underbar p_1^{*}\mtr{\CL}^\vee)\underbar p_2^*\ari{\ch}(\mtr{\CL}^\vee))=\sqrt{\deg(\phi_\CL)}\,\ari{\ch}(\mtr{\CL}^\vee)
$$
which implies that
$$
 p_{2*}(\ari{\ch}(p_1^*\mtr{\CL})\ari{\ch}(\mtr{\CP}))={1\over\sqrt{\deg(\phi_\CL)}}\phi_{\CL,*}\ari{\ch}(\mtr{\CL}^\vee).
$$
Now notice that in the proof of equation \refeq{poincample}, we may assume without restriction of generality
that $\CL$ is relatively generated by its sections, which is to say that the natural morphism $\pi^*\pi_*\CL\to \CL$
is surjective. Indeed, for any $n\geqslant 2$, we have
$$
{1\over g!\sqrt{\deg(\phi_{\CL^{\otimes n}})}}\phi_{\CL^{\otimes n},*}(\ac1(\mtr{\CL}^{\otimes n})^g)=
{1\over g!\sqrt{\deg(\phi_{\CL})}}\phi_{\CL,*}([n]_*\ac1(\mtr{\CL})^g)
$$
and for any $k\geqslant 0$,
\begin{equation}
[n]_*\ac1(\mtr{\CL})^k=n^{-2k}[n]_*[n]^*\ac1(\mtr{\CL})^k=n^{2g-2k}\ac1(\mtr{\CL})^k .
\label{dontknoweq}
\end{equation}
Hence
$$
{1\over g!\sqrt{\deg(\phi_{\CL^{\otimes n}})}}\phi_{\CL^{\otimes n},*}(\ac1(\mtr{\CL}^{\otimes n})^g)=
{1\over g!\sqrt{\deg(\phi_\CL)}}\phi_{\CL,*}(\ac1(\mtr{\CL})^g) .
$$
We may thus harmlessly replace $\mtr{\CL}$ by $\mtr{\CL}^{\otimes n}$, where $n$ is  some
large positive integer. In particular, we may assume (and we do) that the morphism $\pi^*\pi_*\CL\to \CL$ is surjective, since
$\CL$ is relatively ample.
Now let $\CE:=\pi^*\pi_*\CL\otimes\CL^\vee$ and let
$$
P_0^\bullet: \dots\to\Lambda^{r}(\CE)\to\Lambda^{r-1}(\CE)\to\dots\to\CE\to\CO\to 0
$$
be the associated Koszul resolution. Let
$$
P^\bullet_1: 0\to\CP\to p_1^*\CE^\vee\otimes\CP\to\dots\to p_1^*\Lambda^{r-1}(\CE)^\vee\otimes\CP\to p_1^*\Lambda^{r}(\CE)^\vee\otimes\CP
\to\cdots
$$
be the complex $\CP\otimes p_1^*(P_0^\bullet)^\vee$.
All the bundles appearing in the complex $P^\bullet_1$ have natural hermitian metrics and we let $\eta_{\bar P_1}$ be the corresponding Bott-Chern class. Notice the equalities
$$
\eta_{\bar P_1}=\ari{\ch}(\Lambda_{-1}(\mtr{\CE}^\vee))\ari{\ch}(\mtr{\CP})=\ari{c}^{\,\text{top}}(\mtr{\CE})\ari{\Td}^{-1}(\mtr{\CE})\ari{\ch}(\mtr{\CP})
$$
in $\ari{\rm CH}^\bullet(\CA\times_S\CA^\vee)$ (see \cite[last paragraph]{BGS-Complex}). Here
$\Lambda_{-1}(\mtr{\CE}^\vee)$ is the formal $\mZ$-linear combination $\sum_{r\geqslant 0}(-1)^r\Lambda^r(\mtr{\CE}^\vee)$.
Since $\rk(\CE)$ may be assumed arbitrarily large (since we may replace $\CL$ by some of its tensor powers), we see that we may assume that
$\eta_{\bar P_1}=0$ in $\ari{\rm CH}^\bullet(\CA\times_S\CA^\vee)_\mQ$.
Thus we may compute
\begin{equation*}
\begin{split}
p_{2*}(\ari{\ch}(\mtr{\CP}))&=
p_{2*}(\ari{\ch}[-p_1^*\Lambda_{-1}(\mtr{\CE}^\vee)+\mtr{\CO}]\ari{\ch}(\mtr{\CP}))\\
&=
-
\sum_{r=1}^{\rk(\CE)}(-1)^r p_{2*}[\ari{\ch}(\Lambda^r(\pi^*\pi_*(\mtr{\CL})^\vee))
\ari{\ch}(p_1^*\mtr{\CL}^{\otimes r})\ari{\ch}(\mtr{\CP})]\\
&=
-
\sum_{r=1}^{\rk(\CE)}(-1)^r p_{2*}[
\ari{\ch}(p_1^*\mtr{\CL}^{\otimes r})\ari{\ch}(\mtr{\CP})]\ari{\ch}(\Lambda^r(\pi^*\pi_*(\mtr{\CL}))^\vee)\\
&=
-\sum_{r=1}^{\rk(\CE)}(-1)^r {1\over\sqrt{\deg(\phi_{\CL^{\otimes r}})}}\phi_{\CL^{\otimes r},*}(
\ari{\ch}(\mtr{\CL}^{\vee,\otimes r}))\ari{\ch}(\Lambda^r(\pi^*\pi_*(\mtr{\CL}))^\vee)\\
&=
-{1\over\sqrt{\deg(\phi_{\CL})}}\phi_{\CL,*}\Big[\sum_{r=1}^{\rk(\CE)}(-1)^r r^{-g}\big( [r]_*(
\ari{\ch}(\mtr{\CL}^{\vee,\otimes r}))\ari{\ch}(\Lambda^r(\pi^*\pi_*(\mtr{\CL}))^\vee)\big)\Big]\\
&=
-{1\over\sqrt{\deg(\phi_{\CL})}}\phi_{\CL,*}\Big[\sum_{r=1}^{\rk(\CE)}(-1)^r r^{-g}\big([\sum_{s\geqslant 0}r^{2g-2s}\ari{\ch}(\mtr{\CL}^{\vee,\otimes r})^{(s)}]\ari{\ch}(\Lambda^r(\pi^*\pi_*(\mtr{\CL}))^\vee)\big)\Big]\\
&=-{1\over\sqrt{\deg(\phi_{\CL})}}\phi_{\CL,*}\Big[\sum_{r=1}^{\rk(\CE)}\sum_{s\geqslant 0}(-1)^r r^{g-s}\ari{\ch}(\mtr{\CL}^{\vee})^{(s)}\ari{\ch}(\Lambda^r(\pi^*\pi_*(\mtr{\CL}))^\vee)\Big].
\end{split}
\end{equation*}

Now notice that the expression
\begin{equation*}
\begin{split}
[n]_*p_{2*}(\ari{\ch}(\mtr{\CP}))&=p_{2*}((\Id\times [n])_*\ari{\ch}(\mtr{\CP}))\\
&=
 p_{2*}((\Id\times [n])_*(\Id\times [n])^*\sum_{k\geqslant 0}n^{-k}\ari{\ch}(\mtr{\CP})^{(k)})=
p_{2*}(\sum_{k\geqslant 0}n^{2g-k}\ari{\ch}(\mtr{\CP})^{(k)})
\end{split}
\end{equation*}
is a Laurent polynomial in $n\geqslant 2$. The equation \refeq{dontknoweq} shows that the expression
$$
[n]_*\Big(-{1\over\sqrt{\deg(\phi_{\CL})}}\phi_{\CL,*}\Big[\sum_{r=1}^{\rk(\CE)}\sum_{s\geqslant 0}(-1)^r r^{g-s}\ari{\ch}(\mtr{\CL}^{\vee})^{(s)}\ari{\ch}(\Lambda^r(\pi^*\pi_*(\mtr{\CL}))^\vee)\Big]\Big)
$$
is also a Laurent polynomial in $n\geqslant 2$. We may thus identify the coefficients of these polynomials. We obtain the following : if $g+k$ is even, then
$$
p_{2*}(\ari{\ch}(\mtr{\CP}))^{(k)}=
-{1\over\sqrt{\deg(\phi_{\CL})}}\phi_{\CL,*}\Big[\ari{\ch}(\mtr{\CL}^{\vee})^{((g+k)/2)}[\sum_{r=1}^{\rk(\CE)}(-1)^r r^{g-(g+k)/2}\ari{\ch}(\Lambda^r(\pi^*\pi_*(\mtr{\CL}))^\vee)]\Big]
$$
and
$$
p_{2*}(\ari{\ch}(\mtr{\CP}))^{(k)}=0
$$
if $g+k$ is odd. Note that we have already proven the stronger fact that $p_{2*}(\ari{\ch}(\mtr{\CP}))^{(k)}=0$ if $k\not = g$. Thus
$$
p_{2*}(\ari{\ch}(\mtr{\CP}))^{(g)}=
-{1\over\sqrt{\deg(\phi_{\CL})}}\phi_{\CL,*}\Big[\ari{\ch}(\mtr{\CL}^{\vee})^{g}[\sum_{r=1}^{\rk(\CE)}(-1)^r \ari{\ch}(\Lambda^r(\pi^*\pi_*(\mtr{\CL}))^\vee)]\Big].
$$
Using furthermore that
 the left-hand side expression in the last equality is of pure degree $g$ in
$\aCH^\bullet(\CA^\vee)_\mQ$, we deduce that
$$
p_{2*}(\ari{\ch}(\mtr{\CP}))^{(g)}=
-{1\over\sqrt{\deg(\phi_{\CL})}}\phi_{\CL,*}\Big[\ari{\ch}(\mtr{\CL}^{\vee})^{g}[\sum_{r=1}^{\rk(\CE)}(-1)^r {\rk(\CE)\choose r}]\Big].
$$
Now notice that  by the binomial formula $\sum_{r=1}^{\rk(\CE)}(-1)^r {\rk(\CE)\choose r}=(1-1)^{\rk(\CE)}-1=-1$. This proves equation \refeq{poincample}.

\subsection{Proof of \ref{mainth1}.3}

Let $Z$ be an analytic cycle of pure codimension $c$ on $\CA(\mC)$. In view of the assumption on $S$,
$\CA(\mC)$ is a finite disjoint union of abelian varieties and so we may (and do) choose a translation invariant Kähler form
on $\CA(\mC)$. A current $\mfg$ on $\CA(\mC)$ of type $(c-1,c-1)$ is said to be
a harmonic Green current for $Z$ (with respect to the Kähler form), if it satisfies the
following properties : $\mfg$ is a Green current for
$Z$, the differential form $\ddc \mfg+\delta_Z$ is  harmonic and $\int_{\CA(\mC)}\mfg\wedge\kappa=0$
for any harmonic form $\kappa$ on $\CA(\mC)$. Notice now that a differential form on $\CA(\mC)$ is harmonic if and
only if it is translation invariant (see for instance \cite[p. 648]{deRham-Elie-Cartan}). Hence the concept of harmonic Green current is independent of
the choice of the translation invariant Kähler form.

The property (a) in Theorem \ref{mainth0} shows that $\mfg_\CA$ is a Green current for $S_0^\vee(\mC)$ and the  property 2 in Theorem \ref{mainth1} shows that
 $\ddc \mfg+\delta_{S_0^\vee(\mC)}$ is harmonic. Let now $\kappa$ be a harmonic
 form of type $(1,1)$ on $\CA(\mC)$. We know that $\kappa$ is $d$-closed and that
 $[n]^*\kappa=n^2\cdot\kappa$ for any $n\geqslant 2$.  We may thus compute
$$
\int_{\CA(\mC)}\mfg_\CA\wedge\kappa=\int_{\CA(\mC)}[n]_*(\mfg_\CA\wedge\kappa)=
n^{-2}\int_{\CA(\mC)}[n]_*(\mfg_\CA\wedge[n]^*\kappa)=n^{-2}\int_{\CA(\mC)}\mfg_\CA\wedge\kappa
$$
and hence $\int_{\CA(\mC)}\mfg_\CA\wedge\kappa=0.$
Thus $\mfg_\CA$ is harmonic.

\subsection{Proof of \ref{mainth1}.4}

In the next section, we shall give an expression for $\mfg_\CA$, which depends only
on $\CA_\mC$ (see the formula \refeq{dmfg0}). This implies the assertion.

\subsection{Proof of \ref{mainth1}.5}

The proof of 1.2.5 is postponed to the end of the proof of Theorem 1.3.1. See the paragraph before subsection \ref{proof132}.

\subsection{Proof of \ref{mainth1}.6}

This is a direct consequence of Theorem \ref{mainth0}.1 and \cite[Cor. 4.7 (i)]{Burgos-Litcanu-Singular} (thanks to J.-I. Burgos for providing this proof).

\subsection{Proof of \ref{mainth1}.7}

We leave the proof that $\BC(\mC)$ is transverse to $S_0^\vee(\mC)$ to the reader. The equation $\BC^*\mfg_\CA=\mfg_{\CA_T}$
follows from formula \refeq{dmfg0} for $\mfg_\CA$, which will be proved in the next section, together
with \cite[Th. 9.11 (ii)]{Burgos-Litcanu-Singular} and the fact that the higher analytic torsion forms of Bismut-K\"ohler are compatible with base-change.

\section{Proof of Theorem \ref{mainth2}}

\subsection{Proof of \ref{mainth2}.1}

This is the most difficult point to prove. We shall construct a class of currents $\mfg^0_\CA$
which naturally restricts to  the degree $(g-1,g-1)$ part of the analytic torsion and we shall prove that $\mfg_\CA^0$ satisfies the
axioms defining $\mfg_\CA$. The arithmetic Riemann-Roch theorem in higher degrees plays a
crucial role here.

\subsubsection{Definition of $\mfg_\CA^0$}

Let
$$
V: 0\to \CO_\CA\to V_0\to\dots\to V_r\to V_{r+1}\to\cdots
$$
be a resolution of $\CO_\CA$ by $\pi_*$-acyclic vector bundles. Dualising, we get a resolution
$$
V^\vee:\dots\to V_{r+1}^\vee\to V_{r}^\vee\to\dots\to V_0^\vee\to\CO_\CA\to 0
$$
of $\CO_\CA$ on the left.  The first hypercohomology
spectral sequence of the complex $V^\vee\otimes\CP$ for the functor $p_{2,*}$ provides an exact sequence
$$
{\cal H}:\dots\to\R^g p_{2*}(V_r^\vee\otimes\CP)\to\R^g p_{2*}(V_{r-1}^\vee\otimes\CP)\to\dots\to\R^g p_{2*}(V_0^\vee\otimes\CP)\to{\epsilon^\vee}_*(\omega^\vee_{A/S})\to 0 .
$$
Now endow the vector bundles $V_r$ with conjugation-invariant hermitian metrics. The line bundle $\mtr{\omega}_{\CA/S}$ is endowed with its $L^2$-metric. This metric does not depend
on the choice of $\lambda$. This follows from the explicit formula for the $L^2$-metric
on Hodge cohomology given in \cite[Lemma 2.7]{Maillot-Rossler-On-the}.
The arithmetic Riemann-Roch \cite{Gillet-Rossler-Soule-An-arithmetic} in higher degrees applied to $\mtr{\CP}$ and $p_2$ is the identity
\begin{multline*}
(-1)^g\ach(\sum_{r\geqslant 0}(-1)^{r}\R^g p_{2*}(\mtr{V}_r^\vee\otimes\mtr{\CP}))-
\sum_{r\geqslant 0}(-1)^{r}T(\lambda,\mtr{V}_r^\vee\otimes\mtr{\CP})+
\int_{p_2}\Td(\mtr{\T}p_2)\ch(\mtr{\CP})\eta_{\bar{V}^\vee}\\
=p_{2,*}(\ari{\Td}(\mtr{\T}p_2)\ach(\mtr{\CP}))-\int_{p_2}\ch(\CP)R(\T p_2)\Td(\T p_2)
\end{multline*}
in $\ari{\rm CH}^\bullet(\CA^\vee)_\mQ$. Here $\eta_{\bar{V}^\vee}$ is the Bott-Chern secondary class of
$\mtr{V}^\vee$, where $\CO_\CA$ has index $0$. We have identified $\lambda$ with
$p_1^*\lambda$.

Notice first that
\begin{equation*}
\begin{split}
[\int_{p_2}\ch(\CP)R(\T p_2)\Td(\T p_2)]^{(g-1)}&=[\int_{p_2}\ch(\CP)R(\T \pi)\Td(\T \pi)]^{(g-1)}\\
&=
[\pi^{\vee,*}(\epsilon^*(R(\T \pi)\Td(\T \pi)))\int_{p_2}\ch(\CP)]^{(g-1)}=0
\end{split}
\end{equation*}
where we used \refeq{poincample}.

Write $\CT(\mtr{\cal H})$ for the homogenous secondary class in the sense of Bismut-Burgos-Litcanu (see \cite[sec. 6]{Burgos-Litcanu-Singular}) of the resolution ${\cal H}$.
 By its very definition, $-\CT(\mtr{\cal H})^{(g-1)}$ is a class of
Green currents for $S_0^\vee(\mC)$ and it is shown in \cite[Th. 10.28]{Burgos-Litcanu-Singular} that
\begin{equation}
[\ach(\sum_{r\geqslant 0}(-1)^{r}\R^g p_{2*}(\mtr{V}_r^\vee\otimes\mtr{\CP}))]^{(g)}=
(S_0^\vee,-\CT(\mtr{\cal H})^{(g-1)})
\label{BLEq}
\end{equation}
in $\aCH^g(\CA^\vee)_\mQ$ and equation \refeq{poincample} shows that
\begin{multline*}
[p_{2,*}(\ari{\Td}(\mtr{\T}p_2)\ach(\mtr{\CP}))]^{(g)}=[p_{2,*}(\ari{\Td}(\mtr{\T}\pi)\ach(\mtr{\CP}))]^{(g)}\\
=[p_{2,*}(\ach(\mtr{\CP}))\pi^{\vee,*}(\epsilon^*(\ari{\Td}(\mtr{\T}\pi)))]^{(g)}=p_{2,*}(\ach(\mtr{\CP}))^{(g)}
\end{multline*}
hence we are led to the equality
\begin{multline}
p_{2,*}(\ach(\mtr{\CP}))^{(g)}=(-1)^g(S_0^\vee,-\CT(\mtr{\cal H})^{(g-1)})-
\sum_{r\geqslant 0}(-1)^{r}T(\lambda,\mtr{V}_r^\vee\otimes\mtr{\CP})^{(g-1)}\\
+[\int_{p_2}\Td(\mtr{\T}p_2)\ch(\mtr{\CP})\eta_{\bar{V}^\vee}]^{(g-1)} .
\label{ffeq}
\end{multline}
This motivates the definition:
\begin{equation}
\boxed{
\mfg^0_\CA:=
-\CT(\mtr{\cal H})^{(g-1)}+
(-1)^{g+1}\sum_{r\geqslant 0}(-1)^{r}T(\lambda,\mtr{V}_r^\vee\otimes\mtr{\CP})^{(g-1)}+
(-1)^g[\int_{p_2}\Td(\mtr{\T}p_2)\ch(\mtr{\CP})\eta_{\bar{V}^\vee}]^{(g-1)}
}
\label{dmfg0}
\end{equation}
\begin{lemma}
The class of currents $\mfg^0_\CA$ does not depend on the resolution
$V$, nor on the metrics on the bundles $V_r$, nor on the translation invariant Kähler form $\lambda$.
\label{resindep}
\end{lemma}
\beginProof
We first prove that the class of currents $\mfg^0_\CA$ does not depend on
$V$ and that it does not depend on the hermitian metrics or on the bundles $V_r$.

Suppose that there is a second resolution $V'$ dominating $V$ :
\begin{diagram}
 &   & 0 &   & 0 &   & &   & 0 &   & 0 &  & \\
  &   & \dTo &   & \dTo &   & &   & \dTo &   & \dTo &  & \\
V': 0 & \rTo & \CO_\CA & \rTo & V'_0 & \rTo & \cdots& \rTo & V'_r & \rTo & V'_{r+1} & \rTo& \cdots\\
 &   & \dTo^{\rm Id} &   & \dTo &   & &   & \dTo &   & \dTo &  & \\
V: 0 & \rTo & \CO_\CA & \rTo & V_0 & \rTo & \cdots& \rTo & V_r & \rTo & V_{r+1} & \rTo& \cdots\\
  &   & \dTo &   & \dTo &   & &   & \dTo &   & \dTo &  & \\
Q: &   & 0  &   \rTo& Q_0 &\rTo   &\cdots &   \rTo& Q_{r} &\rTo   & Q_{r+1} &\rTo  &\cdots \\
  &   &  &   & \dTo &   & &   & \dTo &   & \dTo &  & \\
&   &   &   & 0 &   & &   & 0 &   & 0 &  & \\
\end{diagram}
By assumption the complex $Q$ is exact and we assume that
its objects are $\pi_*$-acyclic.
We endow everything with hermitian metrics.
We shall write $\CH'$ for the exact sequence
$$
{\cal H}':\dots\to\R^g p_{2*}(V_r^{',\vee}\otimes\CP)\to\R^g p_{2*}(V_{r-1}^{',\vee}\otimes\CP)\to\dots\to\R^g p_{2*}(V_0^{',\vee}\otimes\CP)\to{\epsilon^\vee}_*(\omega^\vee_{A/S})\to 0.
$$
 In order to emphasize the dependence of $\mfg^0_\CA$ on the resolution
$V$ together with the collection of
hermitian metrics on the $V_r$, we shall write $\mfg^0_{\mtr{V}}:=\mfg^0_{\CA,\mtr{V}}$
instead of $\mfg^0_\CA$. Recall that $\eta_{\mtr{V}^\vee}$ is the Bott-Chern secondary class
of the sequence $\mtr{V}^\vee$, with $\mtr{\CO}_\CA$ sitting at the index $0$. We shall
accordingly write $\eta_{\mtr{V}^{',\vee}}$ for the Bott-Chern secondary class
of the sequence $\mtr{V}^{',\vee}$, with $\mtr{\CO}_\CA$ sitting at the index $0$.

By definition, we have
\begin{equation*}
\begin{split}
(-1)^g(\mfg^0_{\mtr{V}}-\mfg^0_{\mtr{V}'})&=(-1)^{g+1}\CT(\mtr{\cal H})^{(g-1)}-
\sum_{r\geqslant 0}(-1)^{r}T(\lambda,\mtr{V}_r^\vee\otimes\mtr{\CP})^{(g-1)}+
[\int_{p_2}\Td(\mtr{\T}p_2)\ch(\mtr{\CP})\eta_{\bar{V}^\vee}]^{(g-1)}\\
&-(-1)^{g+1}\CT(\mtr{\cal H}')^{(g-1)}+
\sum_{r\geqslant 0}(-1)^{r}T(\lambda,\mtr{V}_r^{',\vee}\otimes\mtr{\CP})^{(g-1)}-
[\int_{p_2}\Td(\mtr{\T}p_2)\ch(\mtr{\CP})\eta_{\bar{V}^{',\vee}}]^{(g-1)}.
\end{split}
\end{equation*}
Let now
$$
C_r:0\to Q_r^\vee\otimes\CP\to V^{\vee}_r\otimes\CP\to V^{',\vee}_r\otimes\CP\to 0
$$
be the natural exact sequence. All the bundles appearing on $C_r$ are endowed with natural
hermitian metrics. By the symmetry formula \cite[Th. 2.7, p. 271]{BGS-Bott-Chern}, we may compute
$$
\CT(\mtr{\cal H})-\CT(\mtr{\cal H}')=\wt{\ch}(\R^g\pi_*({Q}^\vee\otimes\CP))-\sum_{r\geqslant 0}(-1)^r
\wt{\ch}(\R^g\pi_*({C}_r)) .
$$
On the other hand, the anomaly formula \cite[Th. 3.10, p. 670]{Bismut-Koehler-Higher} tells us that
\begin{eqnarray*}
&&\sum_{r\geqslant 0}(-1)^{r}T(\lambda,\mtr{V}_r^{\vee}\otimes\mtr{\CP})-
\sum_{r\geqslant 0}(-1)^{r}T(\lambda,\mtr{V}_r^{',\vee}\otimes\mtr{\CP})\\
&=&
\sum_{r\geqslant 0}(-1)^{r}T(\lambda,\mtr{Q}_r^\vee\otimes\mtr{\CP})+(-1)^g\sum_{r\geqslant 0}(-1)^r\wt{\ch}(\R^g\pi_*({C}_r))-\int_{p_2}\sum_{r\geqslant 0}(-1)^r\Td(\mtr{\T}p_2)\wt{\ch}(\mtr{C}_r)
\end{eqnarray*}
and
$$
\sum_{r\geqslant 0}(-1)^{r}T(\lambda,\mtr{Q}_r^\vee\otimes\mtr{\CP})=
\int_{p_2}\Td(\mtr{\T}p_2)\wt{\ch}(\mtr{Q}^\vee\otimes\mtr{\CP})-
(-1)^g\wt{\ch}(\R^g\pi_*({Q}^\vee\otimes\CP)).
$$
Furthermore the symmetry formula \cite[Th. 2.7, p. 271]{BGS-Bott-Chern} again implies that
\begin{multline*}
\int_{p_2}\Td(\mtr{\T}p_2)\eta_{\bar{V}^{',\vee}}\ch(\mtr{\CP})-
\int_{p_2}\Td(\mtr{\T}p_2)\eta_{\bar{V}^{\vee}}\ch(\mtr{\CP})\\
=
-\int_{p_2}\Td(\mtr{\T}p_2)\wt{\ch}(\mtr{Q}^\vee\otimes\mtr{\CP})+\int_{p_2}\sum_{r\geqslant 0}(-1)^r
\Td(\mtr{\T}p_2)\wt{\ch}(\mtr{C}_r).
\end{multline*}
Putting everything together, we see that
\begin{equation*}
\begin{split}
(-1)^g(\mfg^0_{\mtr{V}}-\mfg^0_{\mtr{V}'})&=
\Big[-(-1)^g\wt{\ch}(\R^g\pi_*({Q}^\vee\otimes\CP))+(-1)^g\sum_{r\geqslant 0}(-1)^r
\wt{\ch}(\R^g\pi_*({C}_r))\\
&-
(-1)^g\sum_{r\geqslant 0}(-1)^r\wt{\ch}(\R^g\pi_*({C}_r))+\int_{p_2}\sum_{r\geqslant 0}(-1)^r\Td(\mtr{\T}p_2)\wt{\ch}(\mtr{C}_r)\\
&-\int_{p_2}\Td(\mtr{\T}p_2)\wt{\ch}(\mtr{Q}^\vee\otimes\mtr{\CP})+
(-1)^g\wt{\ch}(\R^g\pi_*({Q}^\vee\otimes\CP))\\
&+
\int_{p_2}\Td(\mtr{\T}p_2)\wt{\ch}(\mtr{Q}^\vee\otimes\mtr{\CP})-\int_{p_2}\sum_{r\geqslant 0}(-1)^r
\Td(\mtr{\T}p_2)\wt{\ch}(\mtr{C}_r)\Big]^{(g-1)}=0 .
\end{split}
\end{equation*}
Homological algebra tells us that there always exists a resolution dominating simultaneously two
other ones. Furthermore, we might assume that this resolution satisfies the above conditions
of $\pi_*$-acyclicity (see \cite[chap XX, par. 3, proof of Th. 3.5, p. 773]{Lang-Algebra}).
 Hence we have proven that
$\mfg_{\mtr{V}}^0$ does not depend on
$V$ and that it does not depend on the hermitian metrics on the bundles $V_r$.

We shall now prove that $\mfg^0_\CA$ does not depend on the choice of $\lambda$.
So let $\lambda'$
be another Kähler fibration, which is translation invariant on the fibres. To emphasize the dependence of $\mfg^0_\CA$ on $\lambda$, let
us write $\mfg^0_\CA:=\mfg^{0,\lambda}$. Write $\mtr{\CH}^\lambda$
(resp. $\mtr{\CH}^{\lambda'}$) for the
sequence $\CH$ together with the hermitian metrics induced by $\lambda$ (resp. $\lambda'$).

We compute as before
\begin{equation*}
\begin{split}
(-1)^g(\mfg^{0,\lambda}&-\mfg^{0,\lambda'})\\
&=
(-1)^{g+1}\CT(\mtr{\cal H}^\lambda)^{(g-1)}-
\sum_{r\geqslant 0}(-1)^{r}T(\lambda,\mtr{V}_r^\vee\otimes\mtr{\CP})^{(g-1)}+
[\int_{p_2}\Td({\T}p_2,\lambda)\ch(\mtr{\CP})\eta_{\bar{V}^\vee}]^{(g-1)}\\
&-(-1)^{g+1}\CT(\mtr{\cal H}^{\lambda'})^{(g-1)}+
\sum_{r\geqslant 0}(-1)^{r}T(\lambda',\mtr{V}_r^{\vee}\otimes\mtr{\CP})^{(g-1)}-
[\int_{p_2}\Td({\T}p_2,\lambda')\ch(\mtr{\CP})\eta_{\bar{V}^{\vee}}]^{(g-1)} .
\end{split}
\end{equation*}
We compute
$$
\CT(\mtr{\cal H}^{\lambda'})-\CT(\mtr{\cal H}^\lambda)=-\sum_{r\geqslant 0}(-1)^r
\wt{\ch}(\R^g\pi_*({V}_r^\vee\otimes{\CP}),\lambda',\lambda)
$$
and by the anomaly formula \cite[Th. 3.10, p. 670]{Bismut-Koehler-Higher}
\begin{multline*}
\sum_{r\geqslant 0}(-1)^{r}T(\lambda',\mtr{V}_r^{\vee}\otimes\mtr{\CP})-
\sum_{r\geqslant 0}(-1)^{r}T(\lambda,\mtr{V}_r^{\vee}\otimes\mtr{\CP})\\
=(-1)^g\sum_{r\geqslant 0}(-1)^r\wt{\ch}(\R^g\pi_*({V}_r^\vee\otimes{\CP}),\lambda',\lambda)
-\int_{p_2}\sum_{r\geqslant 0}(-1)^{r}\wt{\Td}(\lambda',\lambda)\ch(\mtr{V}_r^\vee\otimes\mtr{\CP}).
\end{multline*}
Now notice that by the Leibniz formula (see \cite[6.2, (7)]{Rossler-Adams} for details)
\begin{eqnarray*}
&&\int_{p_2}\big(\Td({\T}p_2,\lambda)-\Td({\T}p_2,\lambda')\big)\ch(\mtr{\CP})\eta_{\bar{V}^\vee}=
-\int_{p_2}\ddc(\wt{\Td}(\lambda',\lambda))\ch(\mtr{\CP})\eta_{\bar{V}^\vee}
\\&=&
-\int_{p_2}\wt{\Td}(\lambda',\lambda)\ddc(\ch(\mtr{\CP})\eta_{\bar{V}^\vee})
=
\sum_{r\geqslant 0}(-1)^r\int_{p_2}\wt{\Td}(\lambda',\lambda)\ch(\mtr{V}_r^\vee\otimes\mtr{\CP})-
\int_{p_2}\wt{\Td}(\lambda',\lambda)\ch(\mtr{\CP}) .
\end{eqnarray*}
Now, by the projection formula, since $\lambda$ and $\lambda'$ are translation invariant, we have
$$
\int_{p_2}\wt{\Td}(\lambda',\lambda)\ch(\mtr{\CP})=\epsilon^{\vee,*}(\wt{\Td}(\lambda',\lambda))\int_{p_2}\ch(\mtr{\CP})
$$
and by the formula \refeq{poincample} we have $[\epsilon^{\vee,*}(\wt{\Td}(\lambda',\lambda))\int_{p_2}\ch(\mtr{\CP})]^{(g-1)}=0$.
Assembling everything, we get that $\mfg^{0,\lambda}-\mfg^{0,\lambda'}=0$.
\endProof

\subsubsection{End of proof of \ref{mainth2}.1}

We keep the notation of the last subsection. By construction, the class of currents $\mfg_\CA^0$ has the property that
\begin{equation}
\mfg^0_\CA|_{\CA^\vee(\mC)\backslash S_0^\vee(\mC)}=(-1)^{g+1}T(\lambda,\mtr{\CP}^0)^{(g-1)}
\label{befnine}
\end{equation}
in
$\widetilde{A}^{g-1,g-1}(\CA^\vee\backslash S_0^\vee)$
and contemplating equation \refeq{BLEq}, we see that the elements of
$\mfg_\CA^0$ are Green currents for $S_0^\vee$. Looking at the equation \refeq{ffeq} we see
that $\mfg^0_\CA$ satisfies (a) and (b) in
Theorem \ref{mainth0}.

Now we want to prove that
\begin{equation}
\mfg^0_\CA=[n]_*\mfg^0_\CA
\label{mfgeq}
\end{equation}
for all $n\geqslant 2$.

Fix $n\geqslant 2$. We claim that to prove equation \refeq{mfgeq}, we may
assume that $\ker\ [n]_\CA$ is a constant diagonalisable subgroup-scheme of $\CA$.
Indeed, both sides of equation \refeq{mfgeq} depend on
 $\CA^\vee_\mC$ only.
 In proving equation \refeq{mfgeq}, we thus may (and do) replace $R$ by its fraction field
 ${\rm Frac}(R)$.
 We may also assume without restriction of generality
 that $S$ is connected, hence integral (since $S$ is regular by assumption). Now let $S'$ be the normalisation of $S$ in the composite of the
 field extensions of the function field $\kappa(S)$ of $S$, which are defined by the residue fields of the
 $n$-torsion points in the generic fibre $\CA_{\kappa(S)}$ of $\CA$. Then $b:S'\to S$ is finite
 (see \cite[II, 6.3.10]{EGA}) and étale (see \cite[Cor. 20.8, p. 147]{Milne-Abelian}).
 Since $[n]_\CA$ is étale, the group
 scheme of $n$-torsion points on $\CA_{S'}$ is then a constant group scheme. If
 we again replace $S'$ by a finite étale cover, we may assume that $\Gamma(S',\CO_{S'})$
 contains the $n$-th roots of $1$ and the group
 scheme of $n$-torsion points on $\CA_{S'}$ then become diagonalisable (and constant).

 Now, by the projection formula the pull-back morphism
 $b^*:\wt{A}^{g-1,g-1}(\CA^\vee_\mR)\to \wt{A}^{g-1,g-1}(\CA^\vee_{S',\mR})$ is injective.
 Furthermore, since an étale finite morphism is a local isomorphism in the category of
 complex manifolds, we have $b^*([n]_*\mfg_\CA^0)=[n]_*(b^*\mfg_\CA^0)$.
 For the same reason, we also have $b^*\mfg_\CA^0=\mfg_{\CA_{S'}}^0$.

 \begin{lemma}\label{critlem}
 There is an isometric isomorphism
$$
[n]_*\CO_{\CA}\simeq\bigoplus_{\CM\in\CA^\vee[n](S)}\CM
$$
where the left-hand side is endowed with its $L^2$-metric and the direct summands on the right-hand side
are endowed with the metric induced by the Poincaré bundle.
\end{lemma}
\beginProof Let us denote by $\CF_\CA:D^b(\CA)\to D^b(\CA^\vee)$ the Fourier-Mukai transformation. If $X$ is a scheme, the category $D^b(X)$ is a full subcategory of the category $D(X)$ derived from the category of sheaves in $\CO_X$-modules.
Its objects are the complexes with a finite number of non-zero homology sheaves, all of which
are coherent.  The functor $\CF_\CA$ is given by the formula
$$
K^\bullet\mapsto \R p_{2,*}(p_1^*K^\bullet \otimes\CP) .
$$
It can be proven that in our situation, there is a natural isomorphism of functors
$$
\CF_{\CA^\vee}\circ\CF_\CA(\cdot)\simeq\big(([-1]^*(\cdot))\otimes\omega^\vee_\CA\big)[-g]
$$
(here we have identified $\CA$ with $(\CA^\vee)^\vee$)
and
$$
\CF_\CA\circ[n]_*(\cdot)\simeq[n]^{*}\circ\CF_\CA(\cdot).
$$
See \cite{Laumon-Transformation} for this. Now we compute
\begin{multline*}
\CF_{\CA^\vee}\circ\CF_{\CA}([n]_*(\CO_\CA))\simeq
[n]_*(\CO_\CA)\otimes\omega_\CA^\vee[-g]\simeq
\CF_{\CA^\vee}\circ[n]^{*}(\CF_\CA(\CO_\CA))\\
\simeq
\CF_{\CA^\vee}([n]^{*}\big(\epsilon_*^\vee(\epsilon^{\vee,*}(\omega^\vee_{\CA}))[-g]\big) )\simeq
\pi_{\CA,*}({\CP}_\CA|_{\CA\times_S\ker\ [n]_{\CA^\vee}})\otimes\omega^\vee_{\CA}[-g]
\end{multline*}
and we thus obtain a canonical isomorphism
\begin{equation}
[n]_*\CO_{\CA}\simeq\bigoplus_{\CM\in\CA^\vee[n](S)}\CM.
\label{blastis1}
\end{equation}
Now we make a different computation.
Let $G$ be a finite group such that $G_S\simeq\ker\ [n]_\CA$. Let $L:=\Gamma(S,\CO_S)$. If
$\chi:G\to L^*$ is a character of $G$ we let $([n]_*\CO_\CA)_\chi$ be the locally free subbundle of $[n]_*\CO_\CA$ which is the largest
subbundle $\CS$ of $[n]_*\CO$ such that the action of $G$ on $\CS$ is given by multiplication by $\chi$. Since $G_L$ is a diagonalisable group scheme over $L$,
this gives a direct sum decomposition
\begin{equation}
[n]_*\CO_\CA=\oplus_\chi([n]_*\CO_\CA)_\chi.
\label{direqsum}
\end{equation}
Now we use the equivariant form of Bismut's relative curvature formula
(see \cite{Bismut-Ma-Familles}). Let $a\in G$ be any non-zero element. Since the
fixed point scheme of $a$ on $\CA$ is empty, we get
\begin{equation}
\ch_a([n]_*(\mtr{\CO}_A))=0.
\label{eqcurv1}
\end{equation}
Here $\ch_a(\cdot)$ is the equivariant Chern character form associated to the action
of $a$. See \cite{Bismut-Ma-Familles} for details. The non-equivariant relative curvature formula gives
\begin{equation}
\ch([n]_*(\mtr{\CO}_A))=\rk([n]_*\CO_A) .
\label{eqcurv2}
\end{equation}
From the equations \refeq{eqcurv1} and \refeq{eqcurv2} and the fact that $\rk([n]_*\CO_A)=\#G$, we deduce (by finite Fourier theory) that all the $([n]_*\CO_\CA)_\chi$ are of
rank $1$ and that
$\c1(([n]_*\mtr{\CO}_\CA)_\chi)=0$, where $([n]_*\mtr{\CO}_\CA)_\chi$ is endowed with the metric induced by $[n]_*\mtr{\CO}_\CA$.

A completely similar computation using
the relative geometric fixed point formula (see \cite{BFQ}) shows that each line bundle $([n]_*\CO_\CA)_\chi$ is actually a torsion line bundle. Notice also that since $G$ acts by isometries on $\CO_A$,
the direct sum decomposition \refeq{direqsum} is an orthogonal direct sum.

Now notice that since $[n]$ is finite and flat there is an isometric equivariant
isomorphism
$$
\epsilon^*([n]_*\mtr{\CO}_\CA)\simeq\mtr{\CO}_\CA^{G}
$$
where ${\CO}_\CA^{G}$ is endowed with the $G$-action induced by the action of $G$ on itself.

This shows that there exist isometric rigidifications $\epsilon^*(([n]_*\mtr{\CO}_\CA)_\chi)\simeq\mtr{\CO}_S$.

Summing up the discussion of the previous paragraphs, we see that there exists an isometric isomorphism of vector bundles
\begin{equation}
[n]_*\mtr{\CO}_\CA=\oplus_\chi([n]_*\mtr{\CO}_\CA)_\chi
\label{blastis2}
\end{equation}
where the direct sum is orthogonal and where each $([n]_*\mtr{\CO}_\CA)_\chi$ is a torsion line bundle, which is isometrically rigidified and which carries a hermitian metric, whose curvature form vanishes.

Now consider the isomorphism of vector bundles
\begin{equation}
\bigoplus_\chi([n]_*\CO_\CA)_\chi\simeq \bigoplus_{\CM\in\CA^\vee[n](S)}\CM
\label{ischi}
\end{equation}
induced by the isomorphisms \refeq{blastis1} and \refeq{blastis2}.

We shall need the

{\bf Claim.} Let $\CM$ and $\CL$ be torsion line bundles on $\CA$. If there is a non-zero morphism of line bundles  (without rigidification)
$\CL\to\CM$ then $\CL$ and $\CM$ are isomorphic.

This follows from the fact
that the Picard functor of $\CA/S$ is representable by a scheme, which is separated over $S$, together
with the fact that a non-trivial torsion bundle on an abelian variety has  no global sections.
Details are left to the reader.

The claim shows that the isomorphism \refeq{ischi} send each $([n]_*\CO_\CA)_\chi$ into
exactly one $\CM\in\CA^\vee[n](S)$ and that this morphism
is an isomorphism of line bundles. After possibly rescaling
the isomorphism $([n]_*\CO_\CA)_\chi\simeq\CM$ by an element of $\Gamma(S,\CO_S^*)$,
we obtain an isomorphism $([n]_*\CO_\CA)_\chi\simeq\CM$ of rigidified line bundles.
Since both $([n]_*\CO_\CA)_\chi$ and $\CM$ are endowed with
the unique metrics, whose curvature form is translation invariant on the fibres of
$\CA(\mC)/S(\mC)$ and which are compatible with the given rigidifications, this
implies that this is an isometric isomorphism. This completes the proof of the lemma. \endProof

As at the beginning of the proof of assertion 1 in Theorem \ref{mainth1}, we have
 a  diagram
\begin{diagram}
\CA\times_S\CA^\vee &\rTo^{\;[n]\times\Id\;\;} & \CA\times_S \CA^\vee & \rTo & \CA^\vee\\
\dTo^{\Id\times[n]}                          &        &                              &         & \dTo^{[n]}\\
\CA\times_S\CA^\vee &        &\rTo                       &         & \CA^\vee
\end{diagram}
such that
$$
(\Id\times[n])^*\CP\simeq ([n]\times\Id)^*\CP
$$
and such that the outer square is cartesian. Write $q:=p_2\circ ([n]\times\Id)$.
To clarify further computations, write $\R^\bullet
p_{2,*}(V^\vee\otimes\CP)$ for
$\cal H$. We now compute
\begin{eqnarray}
(-1)^g[n]^*\mfg_\CA^0&=&[n]^*\Big[(-1)^{g+1}\CT(\R^\bullet p_{2,*}({V}^\vee\otimes{\CP}))-
\sum_{r\geqslant 0}(-1)^{r}T(\lambda,\mtr{V}_r^\vee\otimes\mtr{\CP})+
\int_{p_2}\Td(\mtr{\T}p_2)\ch(\mtr{\CP})\eta_{\bar{V}^\vee}\Big]\nonumber\\
&=&
(-1)^{g+1}\CT(\R^\bullet q_*({V}^\vee\otimes(\Id\times[n])^*{\CP}))-
\sum_{r\geqslant 0}(-1)^{r}T(\lambda,\mtr{V}_r^\vee\otimes(\Id\times[n])^*\mtr{\CP})\nonumber\\
&&+\int_{q}\Td(\mtr{\T}q)\ch((\Id\times[n])^*\mtr{\CP})\eta_{\bar{V}^\vee}\nonumber\\
&=&
(-1)^{g+1}\CT(\R^\bullet q_*({V}^\vee\otimes([n]\times\Id)^*{\CP}))-
\sum_{r\geqslant 0}(-1)^{r}T(\lambda,\mtr{V}_r^\vee\otimes([n]\times\Id)^*\mtr{\CP})\nonumber\\
&&+\int_{q}\Td(\mtr{\T}q)\ch(([n]\times\Id)^*\mtr{\CP})\eta_{\bar{V}^\vee}\nonumber\label{lrleq}
\end{eqnarray}
In view of Lemma \ref{resindep}, we may replace
$\mtr{V}$ by $[n]^*\mtr{V}_\CA$ and
$\lambda$ by $[n]^*\lambda_\CA$ in the string of equalities \refeq{lrleq} without
changing its truth-value. Thus
\begin{equation*}
\begin{split}
(-1)^g[n]^*\mfg_\CA^0&=(-1)^{g+1}\CT(\R^\bullet q_*(([n]\times\Id)^*({V}^\vee\otimes{\CP}))) -
\sum_{r\geqslant 0}(-1)^{r}T([n]^*\lambda,([n]\times\Id)^*(\mtr{V}_{r}^\vee\otimes\mtr{\CP}))\\
&\qquad +\int_{q}\Td(\mtr{\T}q)([n]\times\Id)^*(\ch(\mtr{\CP})\eta_{\bar{V}^\vee})\\
&\stackrel{(*)}{=}
(-1)^{g+1}\CT(\R^\bullet p_{2,*}(([n]\times\Id)_*({\CO})\otimes{V}^\vee\otimes{\CP})) -
\sum_{r\geqslant 0}(-1)^{r}T(([n]\times\Id)_*(\mtr{\CO})\otimes\mtr{V}_{r}^\vee\otimes\mtr{\CP})\\
&\qquad +\int_{p_2}\Td(\mtr{\T}p_2)([n]\times\Id)_*(1)\ch(\mtr{\CP})\eta_{\bar{V}^\vee}\\
&\stackrel{(**)}{=}
\sum_{\tau\in\CA^\vee(S)}\Big[(-1)^{g+1}\CT(\R^\bullet p_{2,*}((\Id\times(\tau\circ\pi))^*(\CP)\otimes{V}^\vee\otimes{\CP}))\\
&\qquad -
\sum_{r\geqslant 0}(-1)^{r}T((\Id\times(\tau\circ\pi))^*(\mtr{\CP})\otimes\mtr{V}_{r}^\vee\otimes\mtr{\CP})\\
&\qquad \qquad +\int_{p_2}\Td(\mtr{\T}p_2)\ch((\Id\times(\tau\circ\pi))^*(\mtr{\CP})\otimes\mtr{\CP})\eta_{\bar{V}^\vee}\Big]\\
&=
\sum_{\tau\in\CA^\vee(S)}\Big[(-1)^{g+1}\tau_{*}(\CT(\R^\bullet p_{2,*}({V}^\vee\otimes\mtr{\CP})))-
\tau_{*}[\sum_{r\geqslant 0}(-1)^{r}T(\mtr{V}_{r}^\vee\otimes\mtr{\CP})]\\
&\qquad +\tau_{*}[\int_{p_2}\Td(\mtr{\T}p_2)(\mtr{\CP})\eta_{\bar{V}^\vee}]\Big]\\
&=
(-1)^g\sum_{\tau\in\CA^\vee(S)}\tau_{*}\mfg^0_{\CA} .
\end{split}
\end{equation*}
The equality (*) is justified by the following
\begin{prop}
The equality
$$
T([n]^*\lambda,([n]\times\Id)^*(\mtr{V}_{r}^\vee\otimes\mtr{\CP}))=
T(([n]\times\Id)_*(\mtr{\CO})\otimes\mtr{V}_{r}^\vee\otimes\mtr{\CP})
$$
is verified for any $r$.
\label{maprop}
\end{prop}
\beginProof (of Proposition \ref{maprop}). This is a direct consequence
of \cite[Intro., Th. 0.1]{Ma-Formes}.
\endProof
For the equality (**), we used Lemma \ref{critlem}.
We may now compute
$$
[n]_*[n]^*\mfg_\CA^0=n^{2g}\mfg_\CA^0=[n]_*(\sum_{\tau\in\CA^\vee(S)}\tau_{*}\mfg^0_{\CA})=
n^{2g}[n]_*\mfg_\CA^0
$$
i.e.
$$
\mfg_\CA^0=[n]_*\mfg_\CA^0.
$$
We have thus proven that $\mfg_\CA^0$ satisfies the conditions (a), (b), (c) in Theorem \ref{mainth0}.
Thus $\mfg_\CA^0=\mfg_\CA$ and looking at equation \refeq{befnine}, we see that we have
almost concluded the proof of Theorem \ref{mainth2}.1. To finish, we quote
\cite{Koehler-Complex}, where it is shown that $T(\lambda,\mtr{\CP}^0)=
\Td^{-1}(\epsilon^*\mtr{\Omega})\gamma$ for some real differential form $\gamma$ of
type $(g-1,g-1)$ on $\CA^\vee\backslash S_0^\vee$. From the above, we see that $\gamma=(-1)^{g+1}\mfg_\CA$ and
we are done.

{\it Proof of Theorem \ref{mainth1}.5}.

We revert to the hypotheses of the introduction and of Theorem \ref{mainth1}.5 (in particular, we do not suppose anymore that
for some $n\geqslant 2$, the group scheme $\ker\,[n]$ is the constant group scheme).

To verify the equation $\mfg_{\CA\times_S\CB}=\pi_{\CA^\vee}^*(\mfg_{\CA})\ast \pi_{\CB^\vee}^*(\mfg_{\CB})$,
it is sufficient to check that the class of currents $\pi_{\CA^\vee}^*(\mfg_{\CA})\ast \pi_{\CB^\vee}^*(\mfg_{\CB})$ satisfies the axioms (a), (b), (c) in Theorem \ref{mainth0}.

The fact that the elements of $\pi_{\CA^\vee}^*(\mfg_{\CA})\ast \pi_{\CB^\vee}^*(\mfg_{\CB})$ are Green currents for the unit
section of $\CA^\vee\times_S\CB^\vee$ follows immediately from the definitions. This settles (a).

To verify (b), we consider the commutative diagram

\begin{diagram}
& &  & & \CA\times\CB\times\CA^\vee\times\CB^\vee& & & & \\
& &  &\ldTo^{p_{\CA\CA^\vee\CB^\vee}}& &\rdTo^{p_{\CA^\vee\CB\CB^\vee}}& & & \\
\CA\times\CA^\vee&\lTo^{q_{\CA\CA^\vee}}& \CA\times\CA^\vee\times\CB^\vee& &\dTo(2,2)^{p_{\CA^\vee\CB^\vee}}& &\CA^\vee\times\CB\times\CB^\vee &\rTo^{r_{\CB\CB^\vee}}&\CB\times\CB^\vee\\
& & &\rdTo_{q_{\CA^\vee\CB^\vee}}& &\ldTo_{r_{\CA^\vee\CB^\vee}}& & & \\
& & & & \CA^\vee\times\CB^\vee & & & &
\end{diagram}

where the morphisms are the obvious ones. Notice that in this diagram, the square is cartesian.
We compute
\begin{eqnarray*}
&&(-1)^{g_{A\times B}} p_{\CA^\vee\CB^\vee,*}(\ari{\ch}(\mtr{\CP}_{\CA\times\CB}))^{g_{A\times B}}\\
&=&(-1)^{(g_A+g_B)} p_{\CA^\vee\CB^\vee,*}(\ari{\ch}(\mtr{\CP}_{\CA\times\CB}))^{(g_A+g_B)}{=}
(-1)^{(g_A+g_B)} p_{\CA^\vee\CB^\vee,*}(\ari{\ch}(\mtr{\CP}_{\CA\times\CB}))\\
&=&
(-1)^{(g_A+g_B)} p_{\CA^\vee\CB^\vee,*}(\ari{\ch}(p_{\CA\CA^\vee\CB^\vee}^*q_{\CA\CA^\vee}^*\mtr{\CP}_{\CA})\cdot\ari{\ch}(p_{\CA^\vee\CB\CB^\vee}^*r_{\CB\CB^\vee}^*\mtr{\CP}_\CB))\\
&=&
(-1)^{(g_A+g_B)} r_{\CA^\vee\CB^\vee,*}(p_{\CA^\vee\CB\CB^\vee,*}(\ari{\ch}(p_{\CA\CA^\vee\CB^\vee}^*q_{\CA\CA^\vee}^*\mtr{\CP}_{\CA})\cdot\ari{\ch}(p_{\CA^\vee\CB\CB^\vee}^*r_{\CB\CB^\vee}^*\mtr{\CP}_\CB)))\\
&=&
(-1)^{(g_A+g_B)} r_{\CA^\vee\CB^\vee,*}(p_{\CA^\vee\CB\CB^\vee,*}(\ari{\ch}(p_{\CA\CA^\vee\CB^\vee}^*q_{\CA\CA^\vee}^*\mtr{\CP}_{\CA}))\cdot\ari{\ch}(r_{\CB\CB^\vee}^*\mtr{\CP}_\CB))\\
&=&
(-1)^{(g_A+g_B)} r_{\CA^\vee\CB^\vee,*}(r_{\CA^\vee\CB^\vee}^*(q_{\CA^\vee\CB^\vee,*}(\ari{\ch}(q_{\CA\CA^\vee}^*\mtr{\CP}_{\CA})))\cdot\ari{\ch}(r_{\CB\CB^\vee}^*\mtr{\CP}_\CB))\\
&=&(-1)^{(g_A+g_B)} r_{\CA^\vee\CB^\vee,*}(\ari{\ch}(r_{\CB\CB^\vee}^*\mtr{\CP}_\CB))\cdot
q_{\CA^\vee\CB^\vee,*}(\ari{\ch}(q_{\CA\CA^\vee}^*\mtr{\CP}_{\CA}))\\
&=&
(-1)^{g_A}\pi_{\CA^\vee}^*(\ari{\ch}(\mtr{\CP}_{\CA}))\cdot (-1)^{g_B}\pi_{\CB^\vee}^*(\ari{\ch}(\mtr{\CP}_{\CB}))
=\pi_{\CA^\vee}^*((S_0^{\vee,\CA},\mfg_{\CA}))\cdot \pi_{\CB^\vee}^*((S_0^{\vee,\CB},\mfg_{\CB}))\\
&=&
(S_0^{\vee,\CA\times\CB},\pi_{\CA^\vee}^*(\mfg_{\CB})\ast\pi_{\CB^\vee}^*(\mfg_\CB))
\end{eqnarray*}
Here we have used the projection formula repeatedly, as well as the fact that direct images in
arithmetic Chow theory are compatible with base-change. We also used Theorem \ref{mainth1}.2.
(in the second equation before last) and the definition of the intersection product in arithmetic Chow theory (in the last equation). This settles (b).

To verify (c) we revert to the hypothesis made at the beginning of this subsection. In particular, we
suppose that for some $n\geqslant 2$, the scheme $\ker\,[n]$ is the constant group scheme in
$\CA$. As explained at the beginning of this subsection, this does not restrict generality.
First notice that by the definition of the symbols $[n]_*$ and $[n]^*$, we have
$$
[n]^*[n]_*\mfg_\CA=\sum_{\tau\in\CA^\vee[n](S)}\tau^*(\mfg_{\CA})
$$
and thus by Theorem \ref{mainth0}(c),
\begin{equation}
[n]^*\mfg_\CA=\sum_{\tau\in\CA^\vee[n](S)}\tau^*(\mfg_{\CA})
\label{sumtau}
\end{equation}
Of course, similar equations hold for $\mfg_\CB$. Notice also that equation \refeq{sumtau} implies
that $[n]_*\mfg_\CA=\mfg_\CA$, as can be see by applying $[n]_*$ to both sides
of equation \refeq{sumtau} (see the calculation made after Proposition \ref{maprop}). Thus the
equation \refeq{sumtau} is actually equivalent to the equation $[n]_*\mfg_\CA=\mfg_\CA$.

Now we may compute
\begin{eqnarray*}
&&[n]^*_{\CA^\vee\times\CB^\vee}(\pi_{\CA^\vee}^*(\mfg_{\CA})\ast\pi_{\CB^\vee}^*(\mfg_\CB))=
\pi_{\CA^\vee}^*([n]^*_{\CA^\vee}(\mfg_{\CA}))\ast\pi_{\CB^\vee}^*([n]^*_{\CB^\vee}(\mfg_\CB))\\
&=&
\pi_{\CA^\vee}^*(\sum_{\tau\in\CA^\vee[n](S)}\tau^*(\mfg_{\CA}))\ast\pi_{\CB^\vee}^*(\sum_{\tau\in\CB^\vee[n](S)}\tau^*(\mfg_{\CB}))\\
&=&
\sum_{\tau\in\CA^\vee[n](S)\times\CB^\vee[n](S)}\tau^*\big(
\pi_{\CA^\vee}^*(\mfg_{\CA})\ast\pi_{\CB^\vee}^*(\mfg_\CB)
\big)
\end{eqnarray*}
which settles (c).
Here we used in the first line the fact that the $\ast$-product is
naturally compatible with finite étale pull-back.

\subsection{Proof of \ref{mainth2}.2}

\label{proof132}

We shall apply the Adams-Riemann-Roch theorem in Arakelov geometry proven in \cite[Th. 3.6]{Rossler-Adams}. Let
$\mtr{\CM}$ be the rigidified hermitian line bundle on $\CA$ corresponding to $\sigma$.
By assumption, there is an isomorphism $\mtr{\CM}^{\otimes n}\simeq\mtr{\CO}_\CA$ of
rigidified hermitian line bundles.
Let $k,l$ be two positive integers such
that $k=l\ (\mod\ n)$. Let $\mtr{\Omega}:=\mtr{\Omega}_{\CA/S}$.
 The theorem  \cite[Th. 3.6]{Rossler-Adams} implies that the identity
 \begin{eqnarray}
\psi^{k}(\R^\bullet\pi_*\mtr{\CM})-\psi^{k}(T(\lambda,\mtr{\CM}))=
\theta^{k}(\epsilon^*\mtr{\Omega})^{-1}\R^\bullet\pi_*(\mtr{M}^{\otimes k})-
\ch(\theta^{k}(\epsilon^*\mtr{\Omega})^{-1})T(\lambda,\mtr{\CM}^{\otimes k})
\label{arr1}
\end{eqnarray}
holds in $\ari{K}_0(S)[{1\over k}]$
and that the identity
 \begin{eqnarray}
\psi^{l}(\R^\bullet\pi_*\mtr{\CM})-\psi^{l}(T(\lambda,\mtr{\CM}))=
\theta^{l}(\epsilon^*\mtr{\Omega})^{-1}\R^\bullet\pi_*(\mtr{M}^{\otimes l})-
\ch(\theta^{l}(\epsilon^*\mtr{\Omega})^{-1})T(\lambda,\mtr{\CM}^{\otimes l})
\label{arr2}
\end{eqnarray}
holds in $\ari{K}_0(S)[{1\over l}]$. Here the symbols $\psi^*(\cdot)$ refer to the Adams operations
acting on arithmetic $K_0$-theory; see \cite[sec. 2 and before Th. 3.6]{Rossler-Adams} for the exact definition. For the definition of the symbol
$\theta^*(\cdot)$, see \cite[sec. 2]{Rossler-Adams}. In the following computations, we shall need the following properties
of these symbols. Define $\phi^t(\cdot)$ to be the additive operator, which sends
a differential form $\eta$ of type $(r,r)$ to the differential form $t^r\cdot\eta$.
It is proven in \cite[Prop. 4.2]{Rossler-Adams} that $\theta^k(\epsilon^*\mtr{\Omega})$ is a unit in $\ari{K}_0(S)[{1\over k }]$ and we have
\begin{equation}
\ch(\theta^k(\epsilon^*\mtr{\Omega}))=k^{\rk(\Omega)}\Td(\epsilon^*\mtr{\Omega}^\vee)\phi^k(\mtr{\Td}^{-1}(\epsilon^*\mtr{\Omega}^\vee))
\label{toddthet}
\end{equation}
(and similarly for $l$ instead of $k$). See \cite[Lemma 6.11]{Rossler-Adams} for this. Secondly,
if $\eta\in\wt{A}(S_\mR)$, then we have
$$
\psi^k(\eta)=k\cdot\phi^k(\eta)
$$
in $\ari{K}_0(S)$. In view of the fact that $\theta^k(\epsilon^*\mtr{\Omega})$ is a unit in $\ari{K}_0(S)[{1\over k }]$, we get
the equation
$$
\ch(\theta^{k}(\epsilon^*\mtr{\Omega}))\psi^{k}(\R^\bullet\pi_*\mtr{\CM})-\ch(\theta^{k}(\epsilon^*\mtr{\Omega}))\psi^{k}(T(\lambda,\mtr{\CM}))=
\R^\bullet\pi_*(\mtr{M}^{\otimes k})-T(\lambda,\mtr{\CM}^{\otimes k})
$$
in $\ari{K}_0(S)[{1\over k }]$  from equation \refeq{arr1}. Similarly, we get
$$
\ch(\theta^{l}(\epsilon^*\mtr{\Omega}))\psi^{l}(\R^\bullet\pi_*\mtr{\CM})-\ch(\theta^{l}(\epsilon^*\mtr{\Omega}))\psi^{l}(T(\lambda,\mtr{\CM}))=
\R^\bullet\pi_*(\mtr{M}^{\otimes l})-T(\lambda,\mtr{\CM}^{\otimes l})
$$
in $\ari{K}_0(S)[{1\over l }]$. Since $k=l\ (\mod\ n)$, we obtain the equation
\begin{multline*}
\ch(\theta^{k}(\epsilon^*\mtr{\Omega}))\psi^{k}(\R^\bullet\pi_*\mtr{\CM})-\ch(\theta^{k}(\epsilon^*\mtr{\Omega}))\psi^{k}(T(\lambda,\mtr{\CM}))\\
=\ch(\theta^{l}(\epsilon^*\mtr{\Omega}))\psi^{l}(\R^\bullet\pi_*\mtr{\CM})-\ch(\theta^{l}(\epsilon^*\mtr{\Omega}))\psi^{l}(T(\lambda,\mtr{\CM}))
\end{multline*}
 in $\ari{K}_0(S)[{1\over kl }]$. In view of equation \refeq{toddthet} and of the fact that
 $\R^r\pi_*{\CM}=0$ for all $r\geqslant 0$, this translates to the equation
\begin{eqnarray*}
&&k^{g}\Td(\epsilon^*\mtr{\Omega}^\vee)\phi^k(\mtr{\Td}^{-1}(\epsilon^*\mtr{\Omega}^\vee))\psi^{k}(T(\lambda,\mtr{\CM}))=
l^{g}\Td(\epsilon^*\mtr{\Omega}^\vee)\phi^l(\mtr{\Td}^{-1}(\epsilon^*\mtr{\Omega}^\vee))\psi^{l}(T(\lambda,\mtr{\CM}))
\end{eqnarray*}
in $\ari{K}_0(S)[{1\over kl }]$.
Recall that in \cite{Koehler-Complex} it is shown that $T(\lambda,\mtr{\CM})=
\Td^{-1}(\epsilon^*\mtr{\Omega})\gamma$, where $\gamma$ is a real differential form
of type $(g-1,g-1)$ on $S$. So we may rewrite
\begin{multline}
k^{g+1}\Td(\epsilon^*\mtr{\Omega}^\vee)\phi^k(\mtr{\Td}^{-1}(\epsilon^*\mtr{\Omega}^\vee))\phi^{k}(\Td^{-1}(\epsilon^*\mtr{\Omega}))\phi^k(\gamma)\\
=
l^{g+1}\Td(\epsilon^*\mtr{\Omega}^\vee)\phi^l(\mtr{\Td}^{-1}(\epsilon^*\mtr{\Omega}^\vee))\phi^{l}(\Td^{-1}(\epsilon^*\mtr{\Omega}))\phi^l(\gamma)
\label{aleq}
\end{multline}
Furthermore
\begin{lemma}
We have $\Td^{-1}(\epsilon^*\mtr{\Omega}^\vee\oplus\epsilon^*\mtr{\Omega})=1$.
\label{lemtodd}
\end{lemma}
\beginProof
We may (and do) assume that $R=\mC$.
Consider the relative Hodge extension
\begin{equation}
0\to \R^0\pi_*(\Omega)\to H^1_\dR(\CA/S)\to R^1\pi_*(\CO_\CA)\to 0
\label{hodex}
\end{equation}
where $H^1_\dR(\CA/S):=R^1\pi_*(\Omega^\bullet_{\CA/S})$ is the first
relative de Rham cohomology sheaf. The sequence \refeq{hodex} is the expression of
the filtration on $R^1\pi_*(\Omega^\bullet_{\CA/S})$, which comes from the relative Hodge
to de Rham spectral sequence. This spectral sequence is the first hypercohomology spectral
sequence of the relative de Rham complex  $\Omega^\bullet_{\CA/S}$ and it degenerates by
\cite[Prop. 5.3]{Deligne-Theoremes}.
The relative form of the GAGA theorem shows that
there is an isomorphism of holomorphic vector bundles \mbox{$H^1_\dR(\CA/S)(\mC)\simeq
(R^1\pi(\mC)_*\mC)\otimes_\mC\CO_{S(\mC)}$} (see \cite[p. 31]{DMOS}) and via this
isomorphism we endow $H^1_\dR(\CA/S)(\mC)$ with the
fibrewise Hodge metric, whose formula is given in \cite[before Lemma 2.7]{Maillot-Rossler-On-the}. Since the metric on $H^1_\dR(\CA/S)(\mC)$ is locally constant by construction, the curvature matrix of the hermitian vector bundle \mbox{$H^1_\dR(\CA/S)(\mC)\simeq
(R^1\pi(\mC)_*\mC)\otimes_\mC\CO_{S(\mC)}$} vanishes. Now the formula in \cite[Lemma 2.7]{Maillot-Rossler-On-the} shows
that in the sequence \refeq{hodex}, the $L^2$-metric on the first term corresponds to
the induced metric and the $L^2$-metric on the end term corresponds to the quotient metric.
We now view the sequence \refeq{hodex} as a sequence of hermitian vector bundles
with the metrics described above. Using \cite[Th. 3.4.1]{Maillot-Un-Calcul}, we see that the secondary class $\widetilde{\Td}$ of the sequence
\refeq{hodex} is $\ddc$-closed. Thus
$$
\Td(R^1\pi_*(\CO_\CA,L^2)\oplus\R^0\pi_*(\Omega,L^2))=\Td((H^1_\dR(\CA/S),{\rm Hodge\ metric}))=1.
$$
Now notice that by relative Lefschetz duality for Hodge cohomology (see \cite[Lemme 6.2]{Deligne-Theoremes}) and
Grothendieck duality, there is
an isomorphism of $\CO_S$-modules $\phi_\lambda:R^1\pi_*(\CO_\CA)\xrightarrow{\sim} \R^0\pi_*(\Omega)^\vee$, which is dependent on
$\lambda$. To describe it, let $S=\Spec \mC$. Under the Hodge-de Rham splitting of the
sequence \refeq{hodex}, the morphism $\phi_\lambda$ is given by the
formula
$$
\omega\mapsto\omega\wedge\lambda^{g-1}\mapsto\int_{\CA}(\omega\wedge\lambda^{g-1})\wedge(\cdot)
$$
(notice that this formula does actually not depend on the splitting).
Now comparing the last formula with the formula for the Hodge metric in \cite[before Lemma 2.7]{Maillot-Rossler-On-the}, we see that, up to a constant factor,
$\phi_\lambda$ induces an isometry between $R^1\pi_*(\CO_\CA,L^2)$
and the dual of the hermitian vector bundle $\R^0\pi_*(\Omega,L^2)$.  To complete the proof, notice that since the volume of the fibres of $\pi(\mC)$ is locally constant (by the assumption on
$\lambda$), the natural isomorphism of vector bundles $\R^0\pi_*(\Omega,L^2)\simeq\epsilon^*\mtr{\Omega}$
is an isometry up to a locally constant factor. Hence the Chern forms of $\R^0\pi_*(\Omega,L^2)$ and $\epsilon^*\mtr{\Omega}$ are the same. This completes the proof.
\endProof

Together with Lemma \ref{lemtodd}, we see that \refeq{aleq} gives
the identity
\begin{eqnarray*}
&&k^{2g}\Td(\mtr{\Omega}^\vee)\gamma=
l^{2g}\Td(\mtr{\Omega}^\vee)\gamma
\end{eqnarray*}
or in other words
\begin{equation}
\boxed{(k^{2g}-l^{2g})T(\lambda,\mtr{\CM})=0}
\label{fincombeq}
\end{equation}
in $\ari{K}_0(S)[{1\over kl }]$.

We shall now show that equation \refeq{fincombeq} implies the identity
\begin{equation}
2g\cdot n\cdot N_{2g}\cdot T(\lambda,\mtr{\CM})=0
\label{rleq}
\end{equation}
in $\ari{K}_0(S)$. For this consider the following combinatorial lemma.


\begin{lemma}
Let $G$ be an abelian group, written additively. Let $c\geqslant
1$ and let $\alpha\in G$. Suppose that for all $k,l>0$ such that $k=l\ (\mod\ n)$, we have
$$
(l^c-k^c)\cdot\alpha=0
$$
in $G[{1\over kl}]$. Then
$$
{\rm order}(\alpha)\ |\ 2\cdot n\cdot c\cdot[\prod_{p\ {\rm prime},\ p\nmid n,\ (p-1) | c} p]
$$
\label{comlem}
\end{lemma}
\beginProof Rephrasing the hypotheses of the lemma, we find that for any  $k, l > 0$ such that $k  = l\ (\mod\ n)$,
there exist integers $a, b \geqslant 0$ (depending on the couple $(k,l)$) such that:
\begin{equation}
\label{blip_blop}
 k^{a} l^{b} (l^{c} - k^{c}) \cdot g = 0.
\end{equation}
For any $(k, l)$ as above, we will denote by $\gamma$ the integer
such that $l = k + \gamma n$. In what follows, $p$ will be a prime
dividing ${\rm order}(\alpha)$ and we will write $\delta_{p}$ (or in short $\delta$
if no confusion can occur) the positive valuation $v_{p}({\rm order}(\alpha))$.

We deduce from equation (\ref{blip_blop}) that:
\[
p^{\delta}\ |\ k^{a}l^{b}(l^{c} - k^{c}).
\]
From now on, $k$ and $l$ will be chosen such that $k, l \not= 0\ (\mod\ p)$,
this implying in particular that the classes of $k$ and $l$ in $\M{Z}/p^{\delta}\M{Z}$  are invertible.
We thus get:
\[
p^{\delta}\ |\  (l^{c} - k^{c})
\]
or equivalently:
\[
(l/k)^{c} = 1\ (\mod\ p^{\delta})
\]
in $(\M{Z}/p^{\delta}\M{Z})^{\ast}$.
We will denote by $C  \subset (\M{Z}/p^{\delta}\M{Z})^{\ast}$
the set of the classes $(l  / k)$ with $l$ and $k$ as above. The restriction
of the map $\varphi : x \mapsto x^{c}$ to $C$
is then identically equal to $1$.  In what follows, in order
to bound $\delta$, we determine the set $C$. We
must distinguish between two cases.
\medskip

Case 1. The prime $p$ doesn't divide $n$.

Taking $k = n$, we have $(l/k) = (k + \gamma n)/k = 1 + \gamma$
for all the integers $\gamma$ satisfying $(1 + \gamma) \wedge p = 1$. We
deduce from this that necessarily:
\[
C = (\M{Z}/p^{\delta}\M{Z})^{\ast}.
\]
We now need the following well-known group isomorphisms,
which we recall for the sake of the exposition (the group laws
are multiplicative on the left-hand side and additive on the right-hand side):
\begin{itemize}
\item
If $p \not= 2$, $(\M{Z}/p^{\delta}\M{Z})^{\ast} \simeq \M{Z}/p^{\delta - 1}(p - 1)\M{Z}$.
\item
If $p = 2$, $(\M{Z}/2\M{Z})^{\ast} \simeq \{0\}$,  $(\M{Z}/4\M{Z})^{\ast} \simeq \M{Z}/2\M{Z}$
and for $\delta \geqslant 3$, $(\M{Z}/2^{\delta}\M{Z})^{\ast} \simeq (\M{Z}/2\M{Z}) \times (\M{Z}/2^{\delta - 2}\M{Z})$,
the projection on the first factor being the obvious reduction map to $(\M{Z}/4\M{Z})^{\ast} \simeq \M{Z}/2\M{Z}$.
\end{itemize}

Let's then discuss the two subcases $p \not= 2$ and $p = 2$ separately.

If $p \not= 2$, the restriction of  $\varphi$ to $C = (\M{Z}/p^{\delta}\M{Z})^{\ast} \simeq \M{Z}/p^{\delta - 1}(p - 1)\M{Z}$
is equal to $1$ identically. On the righten side, $\varphi$ is the multiplication by $\overline{c}$ the class of $c$ in $ \M{Z}/p^{\delta - 1}(p - 1)\M{Z}$.
One must then have $\overline{c} = 0$, and thus:
\[
(p - 1) p^{\delta -1}\ |\ c,
\]
i.e. the two conditions $p - 1 \ | \ c$ and $\delta_{p} \leqslant 1 + v_{p}(c)$.

If $p = 2$, using the same argument as above, the discussion falls into three different subsubcases:
\begin{itemize}
\item $\delta_{2} = 1$, no condition.

\item $\delta_{2} = 2$, we find that $2\ |\ c$.

\item $\delta_{2} \geqslant 3$, we get that $2\ |\ c$ and $2^{\delta_{2} - 2}\ |\ c$,
i.e. $\delta_{2} \leqslant 2 + v_{2}(c)$.
\end{itemize}
Summing up those three subsubcases, we finally conclude that if $c$ is odd then
$\delta_{2} = 1$ and if $c$ is even then $\delta_{2} \leqslant 2 + v_{2}(c)$.
\medskip

Let's now come to the:
\medskip

Case 2. The prime $p$ divide $n$.

Let's define $\beta := v_{p}(n) \geqslant 1$ and $n' := n / p^{\beta}$ and let's
suppose in addition that $\delta > \beta$.

We compute in $(\M{Z}/p^{\beta}\M{Z})^{\ast}$
\[
l / k = (k + \gamma n)/k = 1 + \gamma n / k = 1 + (\gamma n' / k) p^{\beta} = 1
\]
and thus the set $C$ is contained in the kernel $K$ of the
reduction morphism:
\[
(\M{Z}/p^{\delta}\M{Z})^{\ast} \longrightarrow (\M{Z}/p^{\beta}\M{Z})^{\ast}.
\]
The integer $n'$ being prime to $p$, its class in $(\M{Z}/p^{\delta}\M{Z})$ is
invertible and one can take $k = n'$. We find that for any integer $\gamma$:
\[
l / k = 1 + \gamma n / n' = 1 + \gamma p^{\beta}
\]
in  $(\M{Z}/p^{\delta}\M{Z})^{\ast}$, this implying the set equality:
\[
C = K.
\]
We are thus left to determine $K$. Again, we will distinguish between
two subcases:

If $p \not= 2$ we have:
\[
\# K = \# (\M{Z}/p^{\delta}\M{Z})^{\ast} / \# (\M{Z}/p^{\beta}\M{Z})^{\ast} =
p^{\delta -1} (p - 1) / p^{\beta - 1} (p - 1) = p^{\delta - \beta}
\]
from which we deduce first that $K \simeq \M{Z}/p^{\delta - \beta}\M{Z}$
and  then immediately that $p^{\delta - \beta}\  |\  c$, i.e.
\[
\delta \leqslant \beta + v_{p}(c) = v_{p}(n) + v_{p}(c).
\]

If $p = 2$, the discussion now falls into five different subsubcases.
\begin{itemize}
\item $\beta = 1$ and $\delta_{2} = 1$, no condition.
\item $\beta = 1$ and $\delta_{2} = 2$, we find that $2\ |\ c$, i.e. $c$ is even.
\item $\beta = 1$ and $\delta_{2} \geqslant 3$, then we get $2\ |\ c$ and $2^{\delta_{2} - 2}\ |\ c$,
i.e. $\delta_{2} \leqslant 2 + v_{2}(c)$ with $c$ being necessary even.

To summarize those three first cases let's write that when $\beta = 1$,
if $c$ odd then $\delta_{2} = 1$ and if $c$ is even then $\delta_{2} \leqslant 2 + v_{2}(c)$.

\item $\beta = 2$ and $\delta_{2} = 2$, no condition.
\item In all other cases, the kernel $K$ is contained in
$(\M{Z}/2^{\delta_{2} - 2}\M{Z})$ and as a consequence
is cyclic. We thus get:
\[
\# K = \# (\M{Z}/2^{\delta_{2}}\M{Z})^{\ast} / \# (\M{Z}/2^{\beta}\M{Z})^{\ast} =
2^{\delta_{2} -1} / 2^{\beta - 1} = 2^{\delta_{2} - \beta}
\]
and so $2^{\delta_{2} - \beta}\ |\  c$, from which we deduce:
\[
\delta_{2} \leqslant \beta + v_{2}(c) = v_{2}(n) + v_{2}(c).
\]
\end{itemize}
In conclusion of the subcase $p = 2$, we find that:
\[
\delta_{2} \leqslant v_{2}(n) + v_{2}(c) + w_{2},
\]
with  $w_{2} = 1$ if $v_{2}(n) = 1$ and $c$ is even,
and $w_{2} = 0$ otherwise.
\medskip

Putting everything together, we have finally proven that:
\[
{\rm order}(\alpha)\quad | \quad
F_{2}\; \times
\prod_{\textrm{$p$ prime,}\;  p \not= 2,\;  \textrm{$p$ doesn't divide $n$,}\;  (p - 1) | c} p^{1 + v_{p}(c)}
\;\times
\prod_{\textrm{$p$ prime,}\;  p \not= 2,\;  p | n} p^{v_{p}(n) + v_{p}(c)}
\]
where the $F_{2}$ factor is given by the following rules:
\begin{itemize}
\item if $n$ and $c$ are odd then $F_{2} = 2$,
\item if $n$ is odd and $c$ is even then $F_{2} = 2^{2 + v_{2}(c)}$,
\item if $n$ is even then $F_{2} = 2^{v_{2}(n) + v_{2}(c) + w_{2}}$ ;
with $w_{2} = 1$ if $v_{2}(n)  = 1$ and $c$ is even, and $w_{2} = 0$ otherwise.
\end{itemize}

The lemma's statement is then a direct consequence of this (more precise)
assertion.
\endProof

In view of Lemma \ref{comlem}, the set of identities \refeq{fincombeq} implies that
the order of $T(\lambda,\mtr{\CM})$ in $\ari{K}_0(S)$ divides
$2\cdot n\cdot 2g\cdot[\prod_{p\ {\rm prime},\ p\nmid n,\ (p-1) | 2g} \,p]$.

Now use the notations of the last lemma. It is shown in \cite[Appendix B]{Milnor-Stasheff} that the equality
\begin{equation}
2\cdot{\rm denominator}[(-1)^{c+2\over 2}B_c/c]=2\cdot\prod_{p\ {\rm prime},\ (p-1)|c}p^{{\rm ord}_p(c)+1}
\label{vonst}
\end{equation}
holds if $c$ is even. This proves the identity \refeq{rleq} .

To conclude the proof of Theorem \ref{mainth2}.2, recall that there is
an exact sequence
$$
K_1(\CA^\vee)\xrightarrow{-2\reg_\an}\oplus_{p\geqslant 0}\widetilde{A}^{p,p}(\CA^\vee_\mR)
\stackrel{a\;}{\to}\ari{K}_0(\CA^\vee)\to
K_0(\CA^\vee)\to 0
$$
(see \cite[Th. 6.2 (i)]{Gillet-Soule-CharII} for this).

\section{The case of elliptic schemes}
\label{classun}

In this last section, we shall consider elliptic schemes and compare the conclusions of
Theorems \ref{mainth0}, \ref{mainth1} and \ref{mainth2} with
classical results on elliptic units.

So suppose that $\CA$ is of relative dimension $1$, i.e. that $\CA$ is an elliptic scheme over $S$.
Suppose also that the structural morphism $S\to\Spec R$ is the identity on $\Spec R$.
Let $\sigma\in\Sigma$ be an embedding of $R$ into $\mC$. There exists an isomorphism
of complex Lie groups
\begin{equation}
\CA(\mC)_\sigma:=(\CA\times_{R,\sigma}\mC)(\mC)=\mC/(\mZ+\mZ\cdot\tau_\sigma)
\label{idell}
\end{equation}
where $\tau_\sigma\in\mC$ lies in the upper half plane.
\begin{prop} {\rm (a)} The restriction of $\mfg_{\CA^\vee}$ to $\CA(\mC)_\sigma=\mC/(\mZ+\mZ\cdot\tau_\sigma)\backslash\{0\}$ is given by the function
$$
\phi(z)=\phi_{\CA,\sigma}(z):=-2\log|e^{-z\cdot\eta(z)/2}{\rm sigma}(z)\Delta(\tau_{\sigma})^{1\over 12}|
$$
{\rm (b)} Endow $\mC/(\mZ+\mZ\cdot\tau_\sigma)$ with its Haar measure of total
measure $1$. The function $\phi$ then defines an $L^1$-function on $\mC/(\mZ+\mZ\cdot\tau_\sigma)$ and the
restriction of
$\mfg_{\CA^\vee}$ to $\CA(\mC)_\sigma$ is the current $[\phi]$ associated with $\phi$.
\end{prop}
Here $\Delta(\bullet)$ is the discriminant modular form, ${\rm sigma}(z)$ is the Weierstrass sigma-function associated with
the lattice $[1,\tau_{\sigma}]$ and $\eta$ is the quasi-period map associated with
the lattice $[1,\tau_{\sigma}]$, extended $\mR$-linearly to
all of $\mC$ (see \cite[I, Prop. 5.2]{Silverman-Advanced} for the latter).
\beginProof
The formula (a) follows from Theorem \ref{mainth2}.1 and the formula
for the Ray-Singer analytic torsion of a flat line bundle on an elliptic curve given in
\cite[Th. 4.1]{Ray-Singer-Analytic}.

For statement (b), notice that
there exists a Green form $\eta$ of log type along $0$ (see
\cite[Th. 1.3.5, p.106]{Gillet-Soule-Arithmetic}). In this
situation $\eta$ is a real-valued $C^\infty$ function on $\CA(\mC)_\sigma\backslash\{0\}$,
which is locally and hence globally $L^1$ on $\CA(\mC)_\sigma$ (because $\CA(\mC)_\sigma$ is compact).
By definition the current $[\eta]$ associated with $\eta$ is a Green current for $0$ and
by \cite[Lemma 1.2.4]{Gillet-Soule-Arithmetic}, there exists a $C^\infty$ real-valued function $f$
on $\CA(\mC)_\sigma$, such that $\mfg_{\CA^\vee}|_{\CA(\mC)_\sigma}=[\eta]+[f]=[\eta+f]$. Now by construction
the restriction of $\mfg_{\CA^\vee}$ to $\CA(\mC)_\sigma\backslash\{0\}$ is given by the current
associated with the restriction to $\CA(\mC)_\sigma\backslash\{0\}$ of the locally $L^1$-function $\eta+f$.  Since
$\eta+f$ and
$\phi$ are both $C^\infty$ on $\CA(\mC)_\sigma\backslash\{0\}$, they must actually coincide
on $\CA(\mC)_\sigma\backslash\{0\}$.  This proves (b).
\endProof

The distribution relations for $\mfg_{\CA^\vee}$ (i.e. Theorem \ref{mainth0}.3) imply that
the function $\phi(z)$ has the property that
$$
\sum_{w\in \mC/(\mZ+\mZ\cdot\tau_\sigma),\ n\cdot w=z}\phi_{\CA,\sigma}(w)=\phi_{\CA,\sigma}(z).
$$
for any $n\geqslant 2$ and $z\in \mC/(\mZ+\mZ\cdot\tau_\sigma)\backslash\{0\}$. This also follows from the more precise distribution relations given in \cite[Th. 4.1, p. 43]{Kubert-Lang-Modular}.

Now suppose that $R$ is a Dedekind ring. We suppose given
$\underline{z}\in \CA(S)$, whose image is disjoint from the unit section and such that $n\cdot \underline{z}=0$.
Let $\underline{z}_{\sigma}\in \mC/(\mZ+\mZ\cdot\tau_\sigma)$ be the element corresponding to $\underline{z}$.
We compute from the definition that $N_{2}=24$. Theorem
\ref{mainth2}.2 now implies that there exists $u\in R^*$, which does not depend on $\sigma\in\Sigma$, such that
\begin{equation}
\log|\sigma(u)|=24\cdot n\cdot\phi_{\CA,\sigma}(\underline{z}_{\sigma})
\label{receq}
\end{equation}

In particular the real number
$$
\exp(24\cdot n\cdot\phi_{\CA,\sigma}(\underline{z}_{\sigma}))
$$
is an algebraic unit.
If $R$ is the ring of integers of a number field, $n$ has at least two distinct prime
factors and $\Sigma=\{\sigma,\bar\sigma\}$ then \refeq{receq} is also a consequence of \cite[Th. 2.2, p. 37]{Kubert-Lang-Modular}. Notice that \refeq{receq}
overlaps with part of the reciprocity law for elliptic units, if $\CA_{\mtr{{\rm Frac}(R)}}$ is assumed to have complex multiplication by the ring of integers of an imaginary quadratic field
(see \cite[chap. 19, par. 3, Th. 3]{Lang-Elliptic}).
Special instances of elliptic units were first (implicitly) constructed by Eisenstein in his analytic proof of cubic and quartic reciprocity laws (cf. \cite{Eisenstein3} and \cite{Eisenstein1}, \cite{Eisenstein2}). For a thorough historical discussion of Eisenstein's contribution and additional references, see \cite[\S 8]{Lemmermeyer}.

{\bf Remark.} Our proof of the fact that the real number $\exp(24\cdot n\cdot\phi_{\CA,\sigma}(\underline{z}_{\sigma}))
$ is an algebraic unit shows that $24$ naturally comes from a Bernoulli number via
von Staudt's theorem. Indeed, von Staudt's theorem is the main tool in the proof of \refeq{vonst}.
The proof of the fact that the number \mbox{$\exp(24\cdot n\cdot\phi_{\CA,\sigma}(\underline{z}_{\sigma}))
$} (in fact even the number $\exp(12\cdot n\cdot\phi_{\CA,\sigma}(\underline{z}_{\sigma}))
$) is an algebraic unit, which is given in \cite[Th. 2.2, p. 37]{Kubert-Lang-Modular}, does not seem to establish such a link.

\end{document}